\documentclass[fleqn]{mat01}
\usepackage{times,mathtimy,amssymb}%,latexsym}
\newcommand{\dual}[1]{\left\langle#1\right\rangle}
\newtheorem{lem}{Lemma}
\newtheorem{rem}{\it Remark}

\newtheorem{thm}[lem]{\bf Theorem}

\def\newpage{\vfill\eject}
\def\t{\text}
\def\d{\mbox{d}}

\begin{document}

\setcounter{page}{269}
\firstpage{269}

\font\vv=mtsyb at 10pt
\def\bx{\mbox{\vv{\char'242}}}

\title{Non-linear second-order periodic systems with non-smooth potential}

\author{EVGENIA~H~PAPAGEORGIOU and NIKOLAOS S~PAPAGEORGIOU}

\markboth{Evgenia~H~Papageorgiou and Nikolaos S~Papageorgiou}{Periodic systems}

\address{Department of Mathematics, National Technical University, Zografou Campus,
Athens~15780, Greece\\
\noindent E-mail: npapg@math.ntua.gr}

\volume{114}

\mon{August}

\parts{3}

\Date{MS received 13 November 2002; revised 1 October 2003}


\begin{abstract}
In this paper we study second order non-linear periodic systems
driven by the ordinary vector $p$-Laplacian with a non-smooth,
locally Lipschitz potential function. Our approach is variational
and it is based on the non-smooth critical point theory. We prove
existence and multiplicity results under general growth conditions
on the potential function. Then we establish the existence of
non-trivial homoclinic (to zero) solutions. Our theorem appears to
be the first such result (even for smooth problems) for systems
monitored by the $p$-Laplacian. In the last section of the paper
we examine the scalar \hbox{non-linear} and semilinear problem. Our
approach uses a generalized Landesman--Lazer type condition which
generalizes previous ones used in the literature. Also for the
semilinear case the problem is at resonance at any eigenvalue.
\end{abstract}

\maketitle

\section{Introduction}

In a recent paper [28], we proved existence and multiplicity
results for non-linear second-order periodic systems driven by the
one-dimensional $p$-Laplacian and having a non-smooth potential.
Our results there extended to the recent works of Tang [31,32],
who examined semilinear (i.e. $p=2$) systems with smooth
potential. In this paper we continue the study of non-linear,
non-smooth periodic systems. We prove new existence theorems under
more general growth conditions on the non-smooth potential. In
[28] all the results assumed a strict sub-$p$ growth (i.e. strictly
sublinear potential in the semilinear ($p=2$) case). Here the
growth conditions are more general. Also we obtain new
multiplicity results and we also establish the existence of
non-trivial homoclinic solutions. Our approach is variational and
it is based on the non-smooth critical point theory of Chang [4].
Extensions of this theory were obtained recently by
Kourogenis and Papageorgiou [17] and \hbox{Kourogenis} {\it et~al} [18].

Problems with non-differentiable potential which is only locally Lipschitz  in the state
variable $x \in\mathbb{R}^{N}$, are known as `hemivariational inequalities'
and have applications in mechanics and engineering. For details in this direction
we refer to the book of Naniewicz and Panagiotopoulos [27].

In the last decade there has been an increasing interest for problems involving the
one-dimensional $p$-Laplacian or generalizations of it. We refer to the works of
Dang and Oppenheimer [6], Del Pino {\it et~al} [7], Fabry and Fayyad [8], Gasinski and Papageorgiou [9],
Guo [11], Halidias and Papageorgiou [12], Kyritsi {\it et al} [19],
Manasevich and Mawhin [21], Mawhin [22,23] and the references therein.

\section{Mathematical preliminaries}

As we have already mentioned our approach is variational, based on the
non-smooth critical point theory. For the convenience of the
reader, in this section we recall the basic facts from this
theory. It is based on the Clarke subdifferential theory for
locally Lipschitz functions. Let $X$ be a Banach space and
\hbox{$\varphi\!: X\rightarrow\mathbb{R}$}. We say that $\varphi$ is
locally Lipschitz, if for every bounded open set $U\subseteq X$,
there exists a constant $k_{U}>0$ such that
$|\varphi(y)-\varphi(z)| \leq k_{U}||y-z||$ for all $y,z \in U$.
It is a well-known fact from convex analysis that a proper,
convex and lower semicontinuous function \hbox{$\psi\!: X\rightarrow
\overline{\mathbb{R}}= \mathbb{R}\cup \{+\infty\}$} is locally
Lipschitz in the interior of its effective domain $\hbox{dom}\ \psi=\{x
\in X\!\!:\psi(x)< +\infty\}$. In particular an $\mathbb{R}$-valued,
convex and lower semicontinuous function is locally Lipschitz. In
analogy with the directional derivative of a convex function, for
a locally Lipschitz function \hbox{$\varphi\!:\!X\rightarrow\mathbb{R}$}, we
define the generalized directional at derivative $x\in X$ in the
direction $h\in X$, by
\begin{equation*}
\varphi^{0}(x;h)= \mathop {\lim \sup
}\limits_{x'\xrightarrow[{\lambda  \downarrow 0}]{}x} \frac{\varphi
(x'+ \lambda h)-\varphi (x')}{ \lambda}.
\end{equation*}

It is easy to check that the function $h\rightarrow \varphi^{0}(y;h)$ is sublinear, continuous and so by the
Hahn--Banach theorem it is the support function of a non-empty, convex and $w^{*}$-compact set
\begin{equation*}
\partial \varphi (x)=\{ x^{*} \in X^{*}\!\!: ( x^{*} , h ) \leq
\varphi^{0} (x ; h) \quad \t{for} \ \ \t{all} \ \ h \in X \}.
\end{equation*}

The set $\partial\varphi (x)$ is known as the generalized (or Clarke) subdifferential of $\varphi$
at $x\in X$. If $\varphi,\psi\!\!:X\rightarrow \mathbb{R}$ are both locally Lipschitz functions, then
for all $x \in X$ and all $\lambda \in \mathbb{R}$ we have $\partial
(\varphi + \psi)(x) \subseteq
\partial \varphi (x) + \partial \psi (x)$ and $\partial (\lambda
\varphi)(x) = \lambda \partial \varphi (x)$. Moreover, if $\varphi$ is also convex,
then the subdifferential $\partial \varphi$ coincides with the
subdifferential in the sense of convex analysis. Recall that the convex subdifferential
of $\varphi$ is defined by $\partial \varphi (x)=\{ x^{*} \in X^{*}\!: (x^{*} , y-x ) \leq
\varphi (y)-\varphi(x)\ \  \t{for} \ \ \t{all} \ \ y \in X \}$. Also if $\varphi\in C^{1}(X)$,
then $\partial\varphi (x) = \{ \varphi'(x) \}$ for all $x\in X$.

Given a locally Lipschitz function $\varphi\!: X\rightarrow \mathbb{R}$, a point $x\in X$ is said to be a
`critical point' of $\varphi$, if $0\in \partial\varphi(x)$. If $\varphi\in C^{1}(X)$, then as we
saw above, $\partial\varphi (x) = \{ \varphi'(x) \}$ and so this definition of critical point coincides
with the classical (smooth) one. It is easy to see that if $x\in X$ is a local extremum of $\varphi$
(i.e. a local minimum or a local maximum), then $0\in \partial\varphi(x)$. From the smooth critical
point theory, we know that a basic tool is a compactness-type condition, known as the
`Palais--Smale condition' (PS-condition for short). In the present non-smooth setting this condition takes the
following form: `A locally Lipschitz function $\varphi\!\!: X\rightarrow \mathbb{R}$ satisfies the
non-smooth PS-condition, if every sequence $\{x_{n}\}_{n\geq 1}\subseteq X$ such that $\{
\varphi(x_{n})\}_{n\geq 1}$ is bounded and $m (x_{n})= \inf
[\|x^{*}_{n}\|\!\!:x^{*}_{n} \in \partial \varphi (x_{n})] \rightarrow
0$ as $n \rightarrow \infty$, has a strongly convergent
subsequence'. A version of the theory based on a weaker condition known as the `non-smooth
$C$-condition' can be found in Kourogenis and Papageorgiou [17].

A $\lambda \in \mathbb{R}$ is said to be an `eigenvalue' of
minus the $p$-Laplacian with periodic boundary conditions, if the
problem
\begin{equation*}
\left \{  \begin{array}{l} -(||x'(t)||^{p-2}x'(t))'= \lambda ||x(t)||^{p-2}x(t) \quad \hbox{a.e \ on} \ T=[0,b]\\[3pt]
x(0)=x(b),  x'(0)=x'(b), \quad 1< p<\infty
\end{array} \right \},
\end{equation*}
has a non-trivial solution $x\in C^{1}(T,\mathbb{R}^{N})$, known
as corresponding to $\lambda$ `eigenfunction'. Let $S$ denote
the set of these eigenvalues. Evidently $0\in S$ and if
$\lambda  \notin S$, then for every $h\in L^{1}(T,
\mathbb{R}^{N})$ the periodic problem
\begin{equation*}
\left \{  \begin{array}{l} -(||x'(t)||^{p-2}x'(t))'= \lambda ||x(t)||^{p-2}x(t)+h(t) \quad \t{a.e} \ \t{on} \ T=[0,b]\\[3pt]
x(0)=x(b),\quad x'(0)=x'(b)
\end{array} \right \},
\end{equation*}
has at least one solution. Each element of $S$ is non-negative and $0$ is the smallest (first)
eigenvalue. If $N=1$ (scalar case), by direct integration of the equation we obtain all the
eigenvalues which are
\begin{equation*}
\lambda_{n}=\left(\frac{2n\pi_{p}}{b}\right)^{p}, \quad \t{where} \ \ \pi_{p}=2(p-1)^{{1}/{p}}
\frac{({\pi}/{p})}{\sin({\pi}/{p})}.
\end{equation*}

When $p=2$ (semilinear case), then $\pi_{2}=\pi$ and we recover the well-known eigenvalues of the
`scalar periodic negative Laplacian' which are $\lambda_{n}=(n\omega)^{2}$ with
$\omega={2\pi}/{b}$. In the case $N>1$ (vector case),  $\{\lambda_{n}\}_{n\geq 1}\subseteq S$
but $S$ contains more elements (see [22]).

\section{Existence theorem}

In this section we prove an existence theorem for non-smooth
periodic systems driven by the ordinary vector $p$-Laplacian,
which will be used in our investigation of homoclinic orbits in
\S5. It concerns the following non-linear and non-smooth
periodic system:
\begin{equation}
\hskip -4.3pc \left \{  \begin{array}{l}
-(||x'(t)||^{p-2}x'(t))'+g(t)\|x(t)\|^{p-2}x(t)
\in \partial j(t,x(t)) \  \t{a.e} \ \t{on} \ T=[0,b]\\[2pt]
x(0)=x(b)  ,  x'(0)=x'(b),\; 1< p<\infty
\end{array} \right \}.
\end{equation}

Our hypotheses on the data of (1) are the following:\vspace{.4pc}
%\begin{itemize}
%\leftskip .2pc

\noindent $\hbox{H(g)}$: $g\in C(T)$, $g(0)=g(b)$ and for all $t\in T$,
$g(t)\geq c>0$.

\noindent {$\hbox{H}(\hbox{j})_1$:}
$j\!\!:\,T\times\mathbb{R}^{N}\longmapsto \mathbb{R}$
is a functional such that $j(\cdot,0)\in L^{\infty}(T)$, $\int_{0}^{b}j(t,0)\d t\geq 0$ and
\begin{itemize}
\leftskip 2.8pc
\item[(i)] for all $x\in \mathbb{R}^{N}$,
$t\longmapsto j(t,x)$ is measurable;
\item[(ii)] for almost all $ t\in T$, the function $x\longmapsto j(t,x)$
is locally Lipschitz;

\item[(iii)] for almost all
$t\in T$, all $x\in \mathbb{R}^{N}$
and all $u\in \partial j(t,x)$, we have
\begin{equation*}
\hskip -1.2pc \|u\|\le a_{1}(t)+c_{1}(t)\|x\|^{r-1},
\end{equation*}
            $1\le r<+\infty$ with
            $a_{1},c_{1}\in L^{\infty}(T)$;
       \item[(iv)] there exists  $M>0$ such that for almost all
            $t\in T$ and all $x \in \mathbb{R}^{N}$ with
            $\|x\|\geq M$ we have
\begin{equation*}
            \hskip -1.1pc \mu j(t,x)\leq -j^{0}(t,x;-x) \quad \t{with} \ \ \mu >p;
\end{equation*}
        \item[(v)]$\mathop {\lim \sup }\limits_{\|x\| \to \infty }\frac{pj(t,x)}
           {\|x\|^{p}}\leq 0$ uniformly for almost all $t\in T$;
       \item[(vi)] there exists $x_{*}\in \mathbb{R}^{N},\;\|x_*\|\geq M$
            such that $\int_{0}^{b}j(t,x_{*})\d t>0$.
\end{itemize}\vspace{.2pc}

\begin{thm}[\!]
If hypotheses $\hbox{\rm H(g)}${\rm ,} $\hbox{\rm H(j)}_{1}$ hold{\rm ,}
then problem $(1)$ has at least one non-trivial
solution $x\in C^{1}(T,\mathbb{R}^{N})$ with
$\|x'(\cdot)\|^{p-2}x'(\cdot)\in W^{1,r'}(T,\mathbb{R}^{N})$.
\end{thm}

\begin{proof}
Let
$\varphi\!\!:W_{\rm per}^{1,p}(T,\mathbb{R}^{N})\rightarrow\mathbb{R}$ be
the locally Lipschitz function defined by
\begin{equation*}
\varphi(x)=\frac{1}{p}\|x'\|_{p}^{p}+\frac{1}{p}\int_{0}^{b}g(t)\|x(t)\|^{p}\d t
-\int_{0}^{b}j(t,x(t))\d t.
\end{equation*}

First we show that $\varphi$ satisfies the non-smooth
PS-condition. To this end let
 $\{x_{n}\}_{n\geq 1}\subseteq W_{\rm per}^{1,p}(T,\mathbb{R}^{N})$ be a sequence such that
 $|\varphi(x_{n})|\leq M_{1}$ for all $n\geq 1$ and some $M_{1}>0$ and $m(x_{n})\rightarrow 0$.

Since $\partial\varphi(x_n)\subseteq
W_{\rm per}^{1,p}(T,\mathbb{R}^N)^*$ is $w$-compact, the norm
functional in a Banach space is weakly lower semicontinuous and
$W_{\rm per}^{1,p}(T,\mathbb{R}^N)$ is embedded compactly in
$C_{\rm per}(T,\mathbb{R}^N)$, from the Weierstrass theorem we know
that we can find $x_n^*\in\partial\varphi(x_n)$ such that
$m(x_n)=$ \hbox{$\|x_n^*\|, n\geq 1$.} We have
$x_n^*=A(x_n)+g|x_n|^{p-2}x_n-u_n$ with
$A\!\!:W_{\rm per}^{1,p}(T,\mathbb{R}^N)\to
W_{\rm per}^{1,p}(T,\mathbb{R}^N)^*$ being the non-linear operator
defined by
\begin{equation*}
\dual{A(x),y}=\int_0^b\|x'(t)\|^{p-2}(x'(t),y'(t))_{\mathbb{R}^N}\d t,\quad \t{for
all}\;x,y\in W_{\rm per}^{1,p}(T,\mathbb{R}^N)
\end{equation*}

$\left.\right.$\vspace{-1.5pc}

\noindent and $u_n\in L^{r'}(T,\mathbb{R}^N),\;u_n(t)\in\partial
j(t,x_n(t))$ a.e. on $T$ (see [5], pp.~47 and 83). It is
easy to check that $A$ is monotone, demicontinuous; thus maximal
monotone (see [14], p.~309).

Combining hypothesis $\hbox{H}(\hbox{j})_1$(iii) with the Lebourg mean value
theorem (see [20] or p.~41 of [5]), we see that for
almost all $t\in T$ and all $x\in\mathbb{R}^N$,
\begin{equation*}
|j(t,x)|\leq\hat{\alpha}_1(t)+\hat{c}_1(t)\|x\|^r\quad\t{with}\;\;\hat{\alpha}_1,\hat{c}_1\in
L^\infty(T)_+.
\end{equation*}

From the choice of the sequence $\{x_n\}_{n\geq 1}\subseteq
W_{\rm per}^{1,p}(T,\mathbb{R}^N)$, we have
\begin{align*}
&\;\;\;\mu\varphi(x_n)+\dual{x_n^*,-x_n}\leq\;\mu
M_1+\varepsilon_n\|x_n\|\quad \t{with}\;\;\varepsilon_n\downarrow 0\\[5.7pt]
&\;\bigg(\frac{\mu}{p}-1\bigg)\|x'_n\|_p^p+\bigg(\frac{\mu}{p}-1\bigg)\int_0^b
g(t)\|x_n(t)\|^p \d t\\[5.7pt]
&\quad \ \,\,-\int_0^b [(u_n(t),-x_n(t))_{\mathbb{R}^N}+\mu
j(t,x_n(t))]\\[5.7pt]
&\;\;\leq\mu M_1+\varepsilon_n\|x_n\|\\[5.7pt]
\Rightarrow &\;
\bigg(\frac{\mu}{p}-1\bigg)\bigg(\|x'_n\|_p^p+\int_0^b
g(t)\|x_n(t)\|^p \d t \bigg)\\[5.7pt]
&\quad \ \,\,+\int_0^b [-j^0(t,x_n(t);-x_n(t))-\mu
j(t,x_n(t))] \d t \\
&\;\;\leq\mu M_1+\varepsilon_n\|x_n\|.
\end{align*}
\newpage

Using hypotheses $\hbox{H}(\hbox{j})_1$(iii) and (iv), we obtain
\begin{align*}
&\;\;\int_0^b [-j^0(t,x_n(t);-x_n(t))-\mu j(t,x_n(t))] \d t\\
&=\;\int_{\{\|x_n\|<M\}} [-j^0(t,x_n(t);-x_n(t))-\mu j(t,x_n(t))]
\d t\; \\
&\quad\,+ \;\;\int_{\{\|x_n\|\geq M\}} [-j^0(t,x_n(t);-x_n(t))-\mu
j(t,x_n(t))] \d t\\
&\;\geq \;-c_2\quad \t{for some}\;c_2>0\;\;\t{and all }\;n\geq 1.
\end{align*}

Therefore it follows that\vspace{-.2pc}
\begin{align*}
\ & \bigg(\frac{\mu}{p}-1\bigg) \bigg
(\|x_{n}'\|_{p}^{p}+\int_{0}^{b}g(t)\|x_{n}(t)\|^{p}\d t\bigg)
\leq \mu M_{1}+\varepsilon_{n}\|x_{n}\|+c_{2}\\
 &\ \ \,\,\quad \t{for some} \ \ c_{2}>0 \ \ \t{and} \ \ \varepsilon_{n}\downarrow 0,\\
&\Rightarrow \bigg(\frac{\mu}{p}-1\bigg) \bigg(\|x_{n}'\|_{p}^{p} +c\|x_{n}\|_{p}^{p}\bigg)
 \leq c_{3}+\varepsilon_{n}\|x_{n}\| \\
&\ \ \,\quad\,\t{with} \ \ c_{3}=\mu M_{1}+c_{2}>0, \\[4pt]
 &\Rightarrow \|x_{n}\|^{p}\leq c_{4}+\varepsilon'_{n}\|x_{n}\|
\quad \t{for some} \ \ c_{4}>0 \ \ \t{and} \ \
\varepsilon'_{n}\downarrow 0.
\end{align*}

From the last inequality it follows that $\{x_{n}\}_{n\geq
1}\subseteq W_{\rm per}^{1,p}(T,\mathbb{R}^{N})$ is bounded and so by
passing to a subsequence if necessary, we may assume that
$x_n\stackrel{w}\to x$ in $W_{\rm per}^{1,p}(T,\mathbb{R}^N)$ and
$x_n\to x$ in $C_{\rm per}(T,\mathbb{R}^N)$. We have
\begin{align*}
\;|\dual{x_n^*,x_n-x}| &=\;|\dual{A(x_n),x_n-x}\\
&\quad\,-\int_0^b
g(t)\|x_n(t)\|^{p-2}(x_n(t),x_n(t)-x(t))_{\mathbb{R}^N} \d t\\
&\quad\,-\int_0^b (u_n(t),x_n(t)-x(t))_{\mathbb{R}^N}
\d t |\leq\varepsilon_n\|x_n-x\|\\
&\Rightarrow \!\!\,\lim \dual{A(x_n),x_n-x}=0.
\end{align*}

Because $A$ is maximal monotone, it is a generalized pseudomonotone
(see [14], p.~365) and so we have
$\dual{A(x_n),x_n}\to \dual{A(x),x}\Rightarrow
\|x'_n\|_p\to\|x'\|_p$. Because $x'_n\stackrel{w}\to x'$ in
$L^p(T,\mathbb{R}^N)$ and the latter is uniformly convex, from
the Kadec--Klee property (see [14], p.~28), we have
$x'_n\to x'$ in $L^p(T,\mathbb{R}^N)$, hence $x_n\to x$ in
$W^{1,p}_{\rm per}(T,\mathbb{R}^N)$. So $\varphi$ satisfies the
non-smooth PS-condition.

Because of hypothesis $\hbox{H}(\hbox{j})_{1}$(v), given $\varepsilon >0$ we
can find $\delta >0$ such that for almost all $t\in T$ and all
$x\in \mathbb{R}^{N}$ with $\|x\|\leq \delta$ we have $j(t,x)\leq
\frac{\varepsilon}{p}\|x\|^{p}$. On the other hand, hypothesis
H(j)$_{1}$(iii) and the Lebourg mean value theorem, imply that
for almost all $t\in T$ and all $x\in \mathbb{R}^{N}$ with
$\|x\|\geq\delta$ we have $j(t,x)\leq c_{5}\|x\|^{r}$ for some
$c_{5}>0$. So finally for almost all $t\in T$ and all $x\in
\mathbb{R}^{N}$ we can write that $j(t,x)\leq
\frac{\varepsilon}{p}\|x\|^{p}+c_{6}\|x\|^{s}$ for some $c_{6}>0$
and with $s>\max \{r,p\}$. Therefore for every $x\in
W_{\rm per}^{1,p}(T,\mathbb{R}^{N})$ we have
\begin{align*}
\varphi(x)
&=\frac{1}{p}\|x'\|_{p}^{p}+\frac{1}{p}\int_{0}^{b}g(t)\|x(t)\|^{p}\d t
-\int_{0}^{b}j(t,x(t))\d t\\
\ &
\geq\frac{1}{p}\|x'\|_{p}^{p}+\frac{c}{p}\|x\|_{p}^{p}-\frac{\varepsilon}{p}\|x\|_{p}^{p}
-c_{7}\|x\|_{\infty}^{s} \quad \t{for some} \ \ c_{7}>0.
\end{align*}

Because $W_{\rm per}^{1,p}(T,\mathbb{R}^{N})$ is embedded
continuously in $C(T,\mathbb{R}^{N})$, we have
\begin{equation*}
\varphi(x)
\geq\frac{1}{p}(\|x'\|_{p}^{p}+(c-\varepsilon)\|x\|_{p}^p)-c_8\|x\|^s
\quad \t{for some} \ \ c_8>0.
\end{equation*}

Taking $\varepsilon <c$ we obtain that
\begin{equation*}
\varphi(x)\geq c_9\|x\|^p-c_8\|x\|^s \quad \t{for some} \ \ c_9>0.
\end{equation*}

Recall that $s>p$. So we can find $\rho>0$ small so that $\inf
[\varphi(x)\!:\|x\|=\rho]=\xi>0$.

On $\mathbb{R}_{+}\setminus \{0\}$, the function $r\rightarrow
{1}/{r^{\mu}}$ is continuous convex, thus it is locally
Lipschitz. From ([5], p.~48) we have that $r\rightarrow
({1}/{r^{\mu}})j(t,rx)$ is locally Lipschitz on
$\mathbb{R}_{+}\setminus \{0\}$ for almost all $t\in T$
(hypothesis H(j)$_{1}$(ii)) and we have
\begin{equation*}
\partial_{r}\left(\frac{1}{r^{\mu}}j(t,rx)\right)\subseteq -\frac{\mu}{r^{\mu +1}}j(t,rx)+
\frac{1}{r^{\mu}}(\partial_{x}j(t,rx),x)_{\mathbb{R}^{N}}.
\end{equation*}

Using Lebourg's mean value theorem, we can find $\lambda\in
(1,r)$ such that
\begin{align*}
&\frac{1}{r^{\mu}}j(t,rx)\!-\!j(t,x) \in\left(\!\!-\frac{\mu}{\lambda^{\mu
\!+\!1}}j(t,\lambda x)\!+\!\frac{1}{\lambda^{\mu}}(\partial_{x}j(t,\lambda x),x)_{\mathbb{R}^{N}}\right)(r\!-\!1),\\
&\Rightarrow\frac{1}{r^{\mu}}j(t,rx)-j(t,x)
=\frac{r-1}{\lambda^{\mu +1}}(-\mu j(t,\lambda x)
+(\partial_{x}j(t,\lambda x),\lambda x)_{\mathbb{R}^{N}}) \\[.4pc]
&\hskip 8.5pc \geq\frac{r-1}{\lambda^{\mu +1}}(-\mu j(t,\lambda x)-j^{0}
(t,\lambda x;- \lambda x))\geq 0 \\[.4pc]
&\hskip 9pc\, \ \ (\t{see hypothesis} \ \ \hbox{H}(\hbox{j})_{1}(\hbox{iv})) \\[.4pc]
&\Rightarrow r^{\mu}j(t,x)\leq j(t,rx) \quad \t{for almost all} \
\ t\in T,\; \t{all} \ \|x\|\geq M
 \ \ \t{and all} \ \ r\geq 1.
\end{align*}

Choosing $x_*\in \mathbb{R}^N$ as postulated by hypothesis
H(j)$_1$(vi), for $\lambda\geq 1$ large we have
\begin{align*}
\ \varphi(\lambda x_*) &= \int_{0}^{b}g(t)\|\lambda x_*\|^{p}\d t-
\int_{0}^{b}j(t,\lambda x_*)\d t\\[.4pc]
\ &\leq \lambda^p \|g\|_{\infty}\|x_*\|^p b-\lambda^{\mu}
\int_{0}^{b}j(t,x_*)\d t\\[.4pc]
\Rightarrow\  &\ \varphi(\lambda x_*)\rightarrow -\infty \ \ \t{as} \
\ \lambda\rightarrow +\infty \ \ (\t{recall that} \ \ \mu>p).
\end{align*}

Thus we can find $\lambda >0$ large so that $\|\lambda x_*\|>p$
and $\varphi(\lambda x_*)<\xi$. Also note that $\varphi(0)\leq 0$
(recall that $\int_{0}^{b}j(t,0)\d t\geq 0 $). Therefore we can
apply the non-smooth mountain pass theorem (see [4] or
[17]) and obtain $x\in
W_{\rm per}^{1,p}(T,\mathbb{R}^{N})$, $x\neq 0$ such that $0\in
\partial\varphi(x)$.

We have $0\in \partial\varphi(x)\subseteq A(x)-\partial I_{j}(x)$
and so
\begin{equation}
A(x)=u-g\|x\|^{p-2}x \quad \t{with} \  u\in
L^{r'}(T,\mathbb{R}^{N}), u(t)\in
\partial j(t,x(t)) \ \t{a.e. on} \ T.
\end{equation}

Let $\theta\in C_{0}^{\infty}((0,b),\mathbb{R}^{N})$. We have
\begin{align*}
\ &\dual{A(x),\theta}=\int_{0}^{b}(u(t),\theta(t))_{\mathbb{R}^{N}}\d t -\int_0^b g(t)\|x(t)\|
^{p-2}(x(t),\theta(t))_{\mathbb{R}^N} \d t\\
\Rightarrow &
\int_{0}^{b}\|x'(t)\|^{p-2}(x'(t),\theta'(t))_{\mathbb{R}^{N}}\d t=
\int_{0}^{b}(u(t),\theta(t))_{\mathbb{R}^{N}}\d t\\
&\qquad\qquad\qquad\quad\qquad\qquad\quad\qquad -  \int_0^b
g(t)\|x(t)\|^{p-2}(x(t),\theta(t))_{\mathbb{R}^N} \d t.
\end{align*}

Since $(\|x'\|^{p-2}x')'\in W^{-1,q}(T,\mathbb{R}^{N})$ (see [1], p.~50), we have\vspace{-.05cm}
\begin{equation*}
\langle-(\|x'\|^{p-2}x')',\theta\rangle_{0}=\langle u-g\|x\|^{p-2}x,\theta\rangle_{0}
\end{equation*}
with $\dual{\cdot,\cdot}_{0}$ denoting the duality brackets for
the pair
$(W_{0}^{1,p}(T,\mathbb{R}^{N}),W^{-1,q}(T,\mathbb{R}^{N})=W_{0}^{1,p}(T,\mathbb{R}^{N})^{*})$.
Since $C_{0}^{\infty}((0,b),\mathbb{R}^{N})$ is dense in
$W_{0}^{1,p}(T,\mathbb{R}^{N})$, it follows that
\begin{align}
-(\|x'(t)\|^{p-2}x'(t))'+g(t)\|x(t)\|^{p-2}x(t) &=u(t)\quad \t{a.e. on } \ T,\nonumber\\
&\quad\,u\in L^{r'}\!\!(T,\mathbb{R}^{N}).
\end{align}

From (3) it follows that $(\|x'(\cdot)\|^{p-2}x'(\cdot))\in
W^{1,r'}(T,\mathbb{R}^{N})$. Because the map $z\rightarrow
\|z\|^{p-2}z$ is a homeomorphism on $\mathbb{R}^{N}$ onto itself
and $W^{1,r'}(T,\mathbb{R}^{N})\subseteq C(T,\mathbb{R}^{N})$, we
infer that $x'\in C(T,\mathbb{R}^{N})$, hence $x\in
C^{1}(T,\mathbb{R}^{N})$.

Next if $y\in W_{\rm per}^{1,p}(T,\mathbb{R}^{N})$, from Green's inequality (integration by parts),
we have
\begin{align*}
\dual{A(x),y} & =\int_{0}^{b}\|x'(t)\|^{p-2}(x'(t),y'(t))_{\mathbb{R}^{N}}\d t\\
\ & =\|x'(b)\|^{p-2}(x'(b),y(b))_{\mathbb{R}^{N}}-
\|x'(0)\|^{p-2}(x'(0),y(0))_{\mathbb{R}^{N}} \\
\ &\quad\,-\int_{0}^{b}((\|x'(t)\|^{p-2}x'(t))',y(t))_{\mathbb{R}^{N}}\d t.
\end{align*}

Using (2) and (3), we obtain\vspace{-.05cm}
\begin{align*}
\ & \|x'(0)\|^{p-2}(x'(0),y(0))_{\mathbb{R}^{N}} =
\|x'(b)\|^{p-2}(x'(b),y(b))_{\mathbb{R}^{N}}\\[2pt]
&\quad \t{for all $y\in W_{\rm per}^{1,p}(T,\mathbb{R}^{N})$},\\[3pt]
\Rightarrow &\|x'(0)\|^{p-2}x'(0)=\|x'(b)\|^{p-2}x'(b),\\[2pt]
\Rightarrow & x'(0)=x'(b).
\end{align*}

Also since $x\in W_{\rm per}^{1,p}(T,\mathbb{R}^{N})$, $x(0)=x(b)$.
Therefore $x\in C^{1}(T,\mathbb{R}^{N})$ is the desired solution of (1).\hfill QED\!
\end{proof}\vspace{.3pc}

\begin{rem}
{\rm The following function is a
non-smooth potential satisfying hypotheses H(j)$_{2}$ (and does
not satisfy the conditions imposed by Tang [31,32] (for $p=2$)
and
Papageorgiou and Papageorgiou [28]). Again for simplicity we drop the $t$-dependence. We have
\begin{align*}
\hskip -4pc j(x)&=
\begin{cases}
  -\|x\|,& \mbox{if $\|x\|\leq 1$}\\
    \frac{1}{\mu}\|x\|^{\mu}-\|x\|\ln \|x\|+c,& \mbox{if $\|x\|>1$}
 \end{cases}, p<\mu \quad\hbox{and}\quad c=\frac{-\mu-1}{\mu}<0,\\
\hskip -4pc &\Rightarrow -j^{0}(x;-x)=
\begin{cases}
   -\|x\|,& \mbox{if $\|x\|\leq 1$} \\
    \|x\|^{\mu}-\|x\|\ln\|x\|-\|x\|,& \mbox{if $\|x\|>1$}
\end{cases}.
\end{align*}

Since $c<0$, we have $-j^{0}(x,-x)-\mu j(x)\geq (\mu-1)\|x\|\ln
\|x\|-\|x\|\geq  0$ if $\|x\|\geq 1$. Thus hypotheses H(j)$_{1}$
hold.}
\end{rem}

\section{Multiplicity theorems}

In this section we prove a multiplicity result. It concerns an
eigenvalue version of problem~(1).
\begin{equation}
\hskip -2.5pc \left \{  \begin{array}{l} -(||x'(t)||^{p-2}x'(t))'+g(t)\|x(t)\|^{p-2}x(t)
\in \lambda\partial j(t,x(t)) \  \t{a.e} \ \t{on} \ T=[0,b]\\[2pt]
x(0)=x(b)  ,  x'(0)=x'(b), \lambda\in\mathbb{R}.
\end{array} \right \}.
\end{equation}

We prove a multiplicity result for a whole semiaxis of values of the parameter $\lambda\in
\mathbb{R}$. Our hypotheses on the non-smooth potential are the following:
\begin{itemize}
\leftskip 1pc
  \item[{H(j)$_2$:}]
       $j:T\times\mathbb{R}^{N}\longmapsto \mathbb{R}$
       is a functional such that $j(\cdot,0)\in L^{\infty}(T)$ and
%       \begin{itemize}
\hskip 1pc
 \leftskip 1pc
       \item[(i)] for all $x\in \mathbb{R}^{N}$,
               $t\longmapsto j(t,x)$
               is measurable;
       \item[(ii)] for almost all $ t\in T$, the function
               $x\longmapsto j(t,x)$
           is locally Lipschitz;

       \item[(iii)] for almost all
            $t\in T$, all $x\in \mathbb{R}^{N}$
            and all $u\in \partial j(t,x)$, we have
\begin{equation*}
\hskip -3pc \|u\|\le c_{1}(t)(1+\|x\|^{r-1}),
\end{equation*}
            $1\le r<p$ with
            $c_{1}\in L^{\infty}(T)$,
       \item[(iv)] $\int_{0}^{b}j(t,0)\d t=0$ and there exists $x_{0}\in L^{r}(T,\mathbb{R}^{N})$
            such that $\int_{0}^{b}j(t,x_{0}(t))\d t>0$;

        \item[(v)]$\mathop {\lim \sup }\limits_{\|x\| \to \infty }\frac{pj(t,x)}
           {\|x\|^{p}}<0$ uniformly for almost all $t\in T$.
%              \end{itemize}
  \end{itemize}\vspace{.2pc}

\begin{thm}[\!]
If hypotheses $\hbox{\rm H(g)}$ and $\hbox{\rm H(j)}_{2}$ hold{\rm ,}
then there exists $\lambda_{*}>0$ such that
for all $\lambda\geq\lambda_{*}$ problem $(4)$ has at least two
non-trivial solutions $x_{1},x_{2}\in C^{1}(T,\mathbb{R}^{N})$
such that $\|x'_{k}(\cdot)\|^{p-2}x'_{k}(\cdot)\in
W^{1,r'}(T,\mathbb{R}^{N})${\rm ,} $k=1,2$.
\end{thm}

\begin{proof}
For every $\lambda\in \mathbb{R}$ we consider the locally Lipschitz
functional $\varphi_{\lambda}\!\!:W^{1,p}_{\rm per}(T,\mathbb{R}^{N})\rightarrow\mathbb{R}$ defined
by
\begin{equation*}
\varphi_{\lambda}(x)=\frac{1}{p}\|x'\|_{p}^{p}+\frac{1}{p}\int_{0}^{b}g(t)\|x(t)\|^{p}\d t
-\lambda\int_{0}^{b}j(t,x(t))\d t.
\end{equation*}

First we show that $\varphi_{\lambda}$ satisfies the PS-condition. For this purpose, we consider
a sequence $\{x_{n}\}_{n\geq 1}\subseteq W_{\rm per}^{1,p}(T,\mathbb{R}^{N})$ such that
 $|\varphi_{\lambda}(x_{n})|\leq M_{1}$ for all $n\geq 1$ and some $M_{1}>0$ and
$m(x_{n})\rightarrow 0$. As before we can find
$x_{n}^{*}\in\partial\varphi_{\lambda}(x_{n})$ such that
$m(x_{n})=\|x_{n}^{*}\|$ for all $n\geq 1$. For every $n\geq 1$
we have $x_{n}^{*} =A(x_{n})+g\|x_{n}\|^{p-2}x_{n}-\lambda u_{n}$
where $A\!:W_{\rm per}^{1,p}(T,\mathbb{R}^{N}) \rightarrow
W_{\rm per}^{1,p}(T,\mathbb{R}^{N})^{*}$ is as in the proof of Theorem
1 and $u_{n} \in L^{r'}(T,\mathbb{R}^{N})$, $u_{n}(t)\in
\partial j(t,x_{n}(t))$ a.e. on $T$. From hypothesis
H(j)$_{2}$(iii) and the Lebourg mean value theorem, we obtain that
for almost all $t\in T$ and all $x\in \mathbb{R}^{N}$,
$|j(t,x)|\leq c_{2}(t)(1+\|x\|^{r})$ with $c_{2}\in
L^{\infty}(T)$. For every $n\geq 1$, we have
\begin{align*}
\varphi_{\lambda}(x_{n})& =\frac{1}{p}\|x'_{n}\|_{p}^{p}+\frac{1}{p}\int_{0}^{b}g(t)\|x_{n}(t)\|^{p}\d t
-\lambda\int_{0}^{b}j(t,x_{n}(t))\d t\\
\ & \geq\frac{1}{p}\|x'_{n}\|_{p}^{p}+\frac{c}{p}\|x_{n}\|_{p}^{p}-\lambda\|c_{2}\|_{\infty}b
-\lambda c_{3}\|x_{n}\|_{p}^{r},
\end{align*}
for some $c_{3}>0$.\pagebreak

Using Young's inequality with $\varepsilon >0$, we obtain $\lambda c_{3}\|x_{n}\|_{p}^{r}
\leq M_{2}(\varepsilon,\lambda)+\frac{\lambda\varepsilon}{p}\|x_{n}\|_{p}^{p}$ for some
$M_{2}(\varepsilon,\lambda)>0$ (recall that $r<p$). Therefore we obtain
\begin{equation}
\frac{1}{p}\|x_{n}'\|_{p}^{p}+\frac{1}{p}(c-\lambda\varepsilon)\|x_{n}\|_{p}^{p}-c_{4}
(\varepsilon,\lambda)\leq\varphi_{\lambda}(x_{n})\leq M_{1}
\end{equation}
for all $n\geq 1$ and some $c_{4}(\varepsilon,\lambda)>0$.

Choose $\varepsilon >0$ so that $\lambda\varepsilon<c$. From (5)
it follows that $\{x_{n}\}_{n\geq 1}\subseteq
W_{\rm per}^{1,p}(T,\mathbb{R}^{N})$ is bounded. Arguing as in the
last part of the proof of Theorem 1, we conclude that
$\varphi_{\lambda}$ satisfies the non-smooth PS-condition. In fact
from (5) we infer that $\varphi_{\lambda}$ is coercive. Also
exploiting the compact embedding of
$W_{\rm per}^{1,p}(T,\mathbb{R}^{N})$ into $C(T,\mathbb{R}^{N})$
(Sobolev embedding theorem), we can check easily that
$\varphi_{\lambda}$ is sequentially weakly lower semicontinuous.
So from the Weierstrass theorem it follows that there exists
$x_{1}\in W_{\rm per}^{1,p}(T,\mathbb{R}^{N})$ such that
$\varphi_{\lambda}(x_{1})=\inf\varphi_{\lambda}$ and
$0\in\partial\varphi_{\lambda}(x_{1})$.

Next, let $\widehat{\psi}\!\!\!\!:L^{r}(T,\mathbb{R}^{N})\rightarrow
\mathbb{R}$ be the integral functional defined by
$\widehat{\psi}(x)=\int_{0}^{b}j(t,x(t))\d t$. By virtue of
hypothesis H(j)$_{2}$(iv), we have $\widehat{\psi}(x_{0})>0$. The
Sobolev space $W_{\rm per}^{1,p} (T,\mathbb{R}^{N})$ is dense in
$L^{r}(T,\mathbb{R}^{N})$. So we can find \hbox{$y\in W_{\rm per}^{1,p}
(T,\mathbb{R}^{N})$} such that \hbox{$\widehat{\psi}(y)>0$}. Therefore
there exists $\lambda_{*}>0$ large enough such that for all
$\lambda\geq\lambda_{*}$ we have $\varphi_{\lambda}(y) =
\frac{1}{p}\|y'\|_{p}^{p}+\int_{0}^{b}g(t)\|y(t)\|^{p}\d t-\lambda\widehat{\psi}(y)<0$.
Hence
$\varphi_{\lambda}(x_{1})\leq\varphi_{\lambda}(y)<0=\varphi_{\lambda}(0)$,
i.e. $x_{1}\neq 0$. Since $0\in\partial\varphi_{\lambda}(x_{1})$
we verify that $x_{1}\in C^{1}(T,\mathbb{R}^{N})$,
$\|x_{1}'(\cdot)\|^{p-2}x_{1}'(\cdot)\in
W^{1,r'}(T,\mathbb{R}^{N})$ and that it is a non-trivial solution
of~(4).

Because of hypothesis H(j)$_{2}$(v) we can find $\theta >0$ and
$\delta >0$ such that for almost all\break $t\in T$ and all
$\|x\|\leq\delta$, we have $j(t,x)\leq
-\frac{\theta}{p}\|x\|^{p}$. Combining this with the growth
condition on $j$, we obtain that for almost all $t\in T$ and all
$x\in \mathbb{R}^{N}$, $j(t,x)\leq
-\frac{\theta}{p}\|x\|^{p}+c_{5}\|x\|^{s}$ for some $c_{5}>0$ and
with $s>p$. So we can write that
\begin{align*}
\varphi_{\lambda}(x) & \geq\frac{1}{p}\|x'\|_{p}^{p}+\frac{c}{p}\|x\|_{p}^{p}
+\frac{\theta}{p}\|x\|_{p}^{p}-c_{6}\|x\|_{s}^{s} \quad \t{for some} \ \ c_{6}>0\\
\ & \geq c_{7}\|x\|^{p}-c_{8}\|x\|^{s} \quad \t{for some} \ \ c_{7}, c_{8}>0.
\end{align*}

Thus if we choose $0<\rho<\min\{1,\|x_{1}\|\}$ small enough, we can have that
\begin{equation*}
\inf[\varphi_{\lambda}(x)\!\!:\|x\|=\rho]=\gamma>0.
\end{equation*}

Since $\varphi(0)=0$, $x_{1}\neq 0$ and $0<\rho<\|x_{1}\|$, we can apply the non-smooth mountain
pass theorem and obtain $x_{2}\in W_{\rm per}^{1,p}(T,\mathbb{R}^{N})$ such that $0=\varphi(0)<
\gamma\leq \varphi_{\lambda}(x_{2})$, hence $x_{2}\neq 0$, $x_{2}\neq x_{1}$ and $0\in
\partial\varphi_{\lambda}(x_{2})$. As before for $k=1,2$ we can check that $x_{k}\in C^{1}
(T,\mathbb{R}^{N})$, $\|x_{k}(\cdot)\|^{p-2}x_{k}(\cdot)\in
W^{1,r'}(T,\mathbb{R}^{N})$ and it solves (4).\hfill QED
\end{proof}%\vspace{.05pc}

\begin{rem}
{\rm The following non-smooth potential
satisfies hypotheses H(j)$_{2}$
(again we drop the $t$-dependence):
\begin{equation*}
\hskip -3pc j(x)=
\begin{cases}
  -\frac{1}{p}\|x\|^{p},& \mbox{if $\|x\|< 1$} \\[4pt]
   \frac{1}{r}\|x\|^{r}+\cos\|x\|+c,& \mbox{if $\|x\|\geq1$}
\end{cases},\quad r<p \quad \hbox{and}\quad c=\frac{1}{p}-\frac{1}{r}-\cos 1.
\end{equation*}

Note that
\begin{equation*}
\partial j(x)=
\begin{cases}
   -\|x\|^{p-2}x,& \mbox{if $\|x\|< 1$} \\[3pt]
   \t{conv} \{-x,x-(\sin1)x\},& \mbox{if $\|x\|=1$}\\[3pt]
    \|x\|^{r-2}x-\frac{x}{\|x\|}\sin \|x\|, & \mbox{if $\|x\|>1$}.
\end{cases}
\end{equation*}}
\end{rem}

\section{Homoclinic solutions}

In this section we turn our attention to the question of
existence of homoclinic solutions (to 0), for the homoclinic
problem in $\mathbb{R}^{N}$ corresponding to (1). Namely, we
consider the problem:
\begin{equation}
\left \{  \begin{array}{l} -(||x'(t)||^{p-2}x'(t))'+g(t)\|x(t)\|^{p-2}x(t)
\in \partial j(t,x(t)) \quad  \t{a.e} \ \t{on} \ \mathbb{R}\\
\|x(t)\|\rightarrow 0, \|x'(t)\|\rightarrow 0 \ \ \t{as} \ \
|t|\rightarrow\infty,\;1< p <\infty
\end{array} \right \}.
\end{equation}

So far the `homoclinic problem' for second order systems has been studied only in the context
of semilinear equations, primarily with smooth potential. We refer to the works of Grossinho
{\it et~al} [10], Korman and Lazer [16], Rabinowitz [29],Yanheng [34] and the references
therein. Non-smooth semilinear systems were studied only recently by Adly and Goeleven [2] and Hu
[13], using different methods. To our knowledge our result is the first one (even in the
context of smooth systems)  on the existence of homoclinic (to 0) orbits for quasilinear
systems. Our approach is based on that of Rabinowitz [29] (see also [10]).

Our hypotheses on the non-smooth potential are the following:\vspace{.4pc}

{H(j)$_3$:} $j\!:\,\mathbb{R}\times\mathbb{R}^{N}\longmapsto
\mathbb{R}$ is a functional such that $j(t,0)=0$ a.e. on $\mathbb{R}$ and

%\begin{itemize}
%\leftskip 1pc
%  \item
       \begin{itemize}
\leftskip 3.5pc
      \item[(i)] for all $x\in \mathbb{R}^{N}$,
               $t\longmapsto j(t,x)$
               is measurable and $2b$-periodic;
       \item[(ii)] for almost all $ t\in \mathbb{R}$, the function
               $x\longmapsto j(t,x)$
           is locally Lipschitz;

       \item[(iii)] for almost all
            $t\in \mathbb{R}$, all $x\in \mathbb{R}^{N}$
            and all $u\in \partial j(t,x)$, we have
\begin{equation*}
            \|u\|\le a_{1}(t)(1+\|x\|^{p-1}),
\end{equation*}
            with $a_{1}\in L^{\infty}(\mathbb{R})$;
       \item[(iv)] there exists  $M>0$ such that for almost all
            $t\in \mathbb{R}$ and all $x \in \mathbb{R}^{N}$ with
            $\|x\|\geq M$, we have
\begin{equation*}
            \mu j(t,x)\leq -j^{0}(t,x;-x) \ \ \t{with} \ \ \mu >p;
\end{equation*}
        \item[(v)]$\mathop {\lim}\limits_{\|x\| \to 0 }\frac{pj(t,x)}
           {\|x\|^{p}}\leq 0$ uniformly for almost all $t\in \mathbb{R}$;
       \item[(vi)] there exists $x_{0}\in \mathbb{R}^{N}$
            such that $\int_{-b}^{b}j(t,x_{0})\d t>0$.
       \end{itemize}
%  \end{itemize}

\begin{rem}
{\rm Hypothesis H(j)$_{3}$(v) is equivalent
to the following one:
\begin{itemize}
\leftskip 5.2pc
       \item[(v)$'$] $\mathop {\lim}\limits_{\|x\| \to 0}\frac{(u,x)_{\mathbb{R}^{N}}}
           {\|x\|^{p}}\leq 0$ uniformly for almost all $t\in \mathbb{R}$ and all $u\in\partial
           j(t,x)$.
\end{itemize}}
\end{rem}

First we show that $(\hbox{v})\Rightarrow (\hbox{v})'$. From the Lebourg mean value theorem, we know that for
almost all $t\in \mathbb{R}$ and all $x\in\mathbb{R}^{N}\setminus\{0\}$, we have
\begin{align*}
j(t,x)-j\left(t,\frac{x}{2}\right) & =\left(u,\frac{x}{2}\right)_{\mathbb{R}^{N}} \quad \t{with} \ \ u\in \partial
j\left(t,\lambda\frac{x}{2}\right),\\
&\hskip 5.5pc \lambda\in (1,2) \ \ (\t{depending on} \ \ t),\\
\Rightarrow \frac{j(t,x)}{\|x\|^{p}} & =\frac{j(t,\frac{x}{2})}{\|x\|^{p}} +
\frac{(u,\frac{x}{2})_{\mathbb{R}^{N}}}{\|x\|^{p}} \\[.5pc]
\ & =\frac{j(t,\frac{x}{2})}{2^{p}\|\frac{x}{2}\|^{p}}+
\frac{\lambda^{p-1}}{2^{p}}
\frac{(u,\frac{\lambda x}{2})_{\mathbb{R}^{N}}}{
\|\frac{\lambda x}{2}\|^{p}},
\end{align*}\vspace{-2pc}

\begin{equation*}
\hskip 2.6pc\Rightarrow 0  \geq \left(1-\frac{1}{2^{p}}\right)\mathop
{\lim}\limits_{\|x\| \to 0 }\frac{j(t,x)}
           {\|x\|^{p}}\geq \mathop {\lim}\limits_{\|x\| \to 0 }
           \frac{\lambda^{p-1}}{2^{p}}
\frac{(u,\frac{\lambda x}{2})_{\mathbb{R}^{N}}}{
\|\frac{\lambda x}{2}\|^{p}}.
\end{equation*}

As $\|x\|\rightarrow 0$, we have $\lambda\downarrow 1$ and so we
conclude that $\mathop {\lim}\limits_{\|x\| \to
0}\frac{(u,x)_{\mathbb{R}^{N}}}{\|x\|^{p}}\leq 0$ uniformly for
almost all $t\in \mathbb{R}$, i.e. $(\hbox{v})'$ holds.

Next we show that $(\hbox{v})'\Rightarrow (\hbox{v})$. From the previous argument for almost all
$t\in \mathbb{R}$ and all $x \in \mathbb{R}^{N}\setminus\{0\}$ we have
\begin{align*}
&\;\frac{j(t,x)}{\|x\|^{p}}=\frac{j(t,\frac{x}{2})}{2^{p}\|\frac{x}{2}\|^{p}}+
\frac{\lambda^{p-1}}{2^{p}} \frac{(u,\frac{\lambda
x}{2})_{\mathbb{R}^{N}}}{ \|\frac{\lambda x}{2}\|^{p}}\\
\Rightarrow &\;\lim_{\|x\|\to\infty}\frac{j(t,x)}{\|x\|^p}\leq 0
\ \ \left(\t{since} \ 1-\frac{1}{2^p}>0 \ \t{and} \;\;
\lambda\downarrow 1 \ \ \t{as} \;\|x\|\to 0\right),
\end{align*}
and the convergence is uniform for almost all $t\in \mathbb{R}$. So (v) holds. Thus we
have proved that $(\hbox{v})\Leftrightarrow (\hbox{v})'$.

Also the hypothesis on the coefficient function $g$ takes the following form:\\[.3pc]
{H(g)$_{1}$}: $g\in C(\mathbb{R})$, $g$ is 2$b$-periodic and for all $t\in [-b,b]$,
$g(t)\geq c >0$.

\begin{thm}[\!]
If hypotheses $\hbox{\rm H(g)}_{1}$ and $\hbox{\,\rm H(j)}_{3}$ hold{\rm ,}
then there exists a non-trivial homoclinic
solution $x\in C(\mathbb{R},\mathbb{R}^{N}) \cap
W^{1,p}(\mathbb{R},\mathbb{R}^{N})$ for problem $(6)$.
\end{thm}

\begin{proof}
We consider the following auxiliary periodic problem:
\begin{equation}
\left \{  \begin{array}{l} -(||x'(t)||^{p-2}x'(t))'+g(t)\|x(t)\|^{p-2}x(t)
\in \partial j(t,x(t))\\
 \qquad\qquad  \t{a.e} \ \t{on} \ T[-nb,nb],\\[4pt]
x(-nb)=x(nb)  ,  x'(-nb)=x'(nb),\; 1< p<\infty
\end{array} \right \}.
\end{equation}

From Theorem~2, we know that problem (7) has a non-trivial
solution $x_{n}\in W_{\rm per}^{1,p}(T_{n},\mathbb{R}^{N})$.\ Let
$\varphi_{n}\!\!:W_{\rm per}^{1,p}(T_{n},\mathbb{R}^{N})\rightarrow
\mathbb{R}$ be the locally Lipschitz energy functional
corresponding to problem (7), i.e.
\begin{equation*}
\varphi_{n}(x)=\frac{1}{p}\|x'\|_{p}^{p}+\frac{1}{p}\int_{-nb}^{nb}g(t)\|x(t)\|^{p}\d t
-\int_{-nb}^{nb}j(t,x(t))\d t.
\end{equation*}

Hereafter by $L_{n}^{p}$ we shall denote the Lebesgue space
$L^{p}(T_{n},\mathbb{R}^{N})$ and by $W^{1,p}_{n}$ the Sobolev
space $W^{1,p}(T_{n},\mathbb{R}^{N})$.

Consider the integral functional $\psi\!\!:L_{1}^{p}\rightarrow
\mathbb{R}$ defined by $\psi(x)= \int_{-b}^{b}j(t,x(t))\d t$. By
virtue of hypothesis H(j)$_{3}$(vi) we have that $\psi(x_{0})>0$.
Because $\psi$ is continuous and
$W^{1,p}_{0}(T_{1},\mathbb{R}^{N})$ is dense in $L_{1}^{p}$, we
can find $\overline{x}\in W^{1,p}_{0}(T_{1},\mathbb{R}^{N})$ such
that $\psi(\overline{x})= \int_{-b}^{b}j(t,\overline{x}(t))\d t>0$.
Then recalling that for almost all $t\in T$, all $\|x\|\geq 1$ and
all $\lambda \geq 1$, we have $\lambda^\mu j(t,x)\leq j(t,\lambda
x)$ we can easily see that
\begin{align*}
\varphi_{1}(\lambda \bar{x})
&=\frac{\lambda^{p}}{p}\|\bar{x}'\|_{p}^{p}+\frac{\lambda^{p}}{p}
\int_{-b}^{b}g(t)\|\bar{x}(t)\|^{p}\d t\\[3pt]
&\quad\, -\int_{-b}^{b}j(t,\lambda\bar{x}(t))\d t \to - \infty \quad \t{as} \
\lambda\to +\infty.
\end{align*}

(Recall $\mu>p$.) So we can find $\lambda_{0}\geq 1$ such that for
all $\lambda\geq\lambda_{0}$ we have $\varphi_{1}(\lambda
\bar{x})<0$. Define $\hat{x}\in
W^{1,p}_{0}(T_{n},\mathbb{R}^{N})$ as follows:
\begin{equation*}
\hat{x}(t)= \begin{cases}
  \bar{x}(t), & \mbox{if $t\in T_{1}$} \\[3pt]
    0, & \mbox{if $t\in T_{n}\setminus T_{1}$}
\end{cases}.
\end{equation*}
Then we have $\varphi_{n}(\lambda \hat{x})=\varphi_{1}(\lambda \bar{x})<0$ for all $\lambda
\geq \lambda_{0}$ (recall that $j(t,0)=0$ a.e. on $\mathbb{R}$).

From the proof of Theorem~1 we know that the solution $x_{n}\in
W^{1,p}_{\rm per}(T_{n},\mathbb{R}^{N})$ of problem (7) is obtained
via the non-smooth mountain pass theorem and so it satisfies (see [17])
\begin{align*}
c_{n}=\mathop {\inf \sup }\limits_{\gamma  \in \Gamma _n t \in [0,1]}
\varphi_{n}(\gamma(t))&=\varphi_{n}(x_{n})\geq \inf [\varphi_{n}(x)\!\!:\|x\|=\rho_{n}]=\xi_{n}>0\\
&\quad \t{and} \ \ 0\in \partial\varphi_{n}(x_{n}),
\end{align*}
where $\Gamma_{n}=\{\gamma\in C([0,1],W_{n}^{1,p})\!:\gamma(0)=0,\gamma(1)=\lambda\hat{x}\}$ with
$\lambda\geq\lambda_{0}$. By continuous extension by constant, we see that for $n_{1}\leq n_{2}$,
we have
\begin{equation*}
W^{1,p}_{n_{1}}\subseteq W^{1,p}_{n_{2}}, \quad \Gamma_{n_{1}}\subseteq \Gamma_{n_{2}}
\quad \t{and so} \ \ c_{n_{2}}\leq c_{n_{1}}.
\end{equation*}

Therefore we have produced a decreasing sequence $\{c_{n}\}_{n\geq 1}$ of critical values.
For every $n\geq 1$ we have
\begin{align}
c_{n}=\varphi_{n}(x_{n})&=\frac{1}{p}\|x'_{n}\|_{L_{n}^{p}}^{p}+\frac{1}{p}\int_{-nb}^{nb}g(t)
\|x_{n}(t)\|^{p}\d t\nonumber\\
&\quad\, -\int_{-nb}^{nb}j(t,x_{n}(t))\d t \leq c_{1}.
\end{align}

Since $0\in \partial\varphi_{n}(x_{n})$, we can find $x_{n}^{*}\in \partial\varphi_{n}(x_{n})$
such that $x_{n}^{*}=0$. So we have
\begin{align}
&A(x_{n})+g\|x_{n}\|^{p-2}x_{n}=u_{n},\quad \t{with} \ \ u_{n}\in L_{n}^{\infty},\nonumber \\
&\quad u_{n}(t)\in\partial j(t,x_{n}(t)) \quad
\t{a.e. on} \ \ T_{n}.
\end{align}

We take the duality brackets (for the pair
($W^{1,p}_{n},(W^{1,p}_{n})^{*}$)) of (9) with $-x_{n}$. We obtain
\begin{align}
-\|x'_{n}\|_{L_{n}^{p}}^{p}-\frac{1}{p}\int_{-nb}^{nb}g(t)
\|x_{n}(t)\|^{p}\d t &=\int_{-nb}^{nb}(u_{n}(t),-x_{n}(t))_{\mathbb{R}^{N}}\d t \nonumber\\
\ &\leq\int_{-nb}^{nb}j^{0}(t,x_{n}(t);-x_{n}(t))\d t.
\end{align}

Multiply (8) with $\mu >p$ and then add to (10). We obtain
\begin{align}
\ &\left(\frac{\mu}{p}-1\right)\|x'_{n}\|_{L_{n}^{p}}^{p}+\left(\frac{\mu}{p}-1\right)\int_{-nb}^{nb}g(t)
\|x_{n}(t)\|^{p}\d t \nonumber\\
\ &\qquad\qquad +\int_{-nb}^{nb}(-j^{0}(t,x_{n}(t);-x_{n}(t))-\mu
j(t,x_{n}(t))\d t\leq \mu c_{1}.
\end{align}

Using hypothesis H(j)$_{3}$(iv), for every $n\geq 1$ we have
\begin{align*}
\ & \int_{-nb}^{nb}(-j^{0}(t,x_{n}(t);-x_{n}(t))-\mu j(t,x_{n}(t))\d t\\
= & \int_{T_{n}\cap\{\|x_{n}(t)\| <M\}}(-j^{0}(t,x_{n}(t);-x_{n}(t))-\mu j(t,x_{n}(t))\d t\\
\ & +\int_{T_{n}\cap\{\|x_{n}(t)\| \geq M\}}(-j^{0}(t,x_{n}(t);-x_{n}(t))-\mu j(t,x_{n}(t))\d t
\geq -\xi_{1},
\end{align*}
for some $\xi_{1}>0$ independent of $n\geq 1$. Using this lower
bound in (11), we obtain
\begin{equation*}
\left(\frac{\mu}{p}-1\right)\|x_{n}\|_{L_{n}^{p}}^{p}+\left(\frac{\mu}{p}-1\right)\int_{-nb}^{nb}g(t)
\|x_{n}(t)\|^{p}\d t\leq c_{1} + \xi_{1}=\xi_{2}
\end{equation*}
with $\xi_{2}>0$ independent of $n\geq 1$. So it follows that
\begin{equation}
\|x_{n}\|_{W_{n}^{1,p}}\leq \xi_{3},
\end{equation}
with $\xi_{3}>0$ independent of $n\geq 1$. Moreover, as in ([29], p.~36), we can have that
\begin{equation}
\|x_{n}\|_{L_{n}^{\infty}}\leq \xi_{4},
\end{equation}
with $\xi_{4}>0$ independent of $n\geq 1$. We extend by
periodicity $x_{n}$ and $u_{n}$ to all of $\mathbb{R}$. From\break (12)
and since $W_{n}^{1,p}$ is embedded compactly in $C_{n}=
C(T_{n},\mathbb{R}^{N})$, by passing to a subsequence  if
necessary, we may assume that $x_{n}\rightarrow x$ in
$C_{\t{loc}}(\mathbb{R},\mathbb{R}^{N})$, hence $x\in C(\mathbb{R}
,\mathbb{R}^{N})$. Also because of hypothesis H(j)$_{3}$(iii), we
have
\begin{align*}
\|u_{n}(t)\|\leq\|a_{1}\|_{\infty}(1+\|x_{n}(t)\|^{r-1})&\leq
\|a_{1}\|_{\infty}(1+\xi_{4}^{r-1}) =\xi_{5}\\
&\quad \t{a.e. on} \ \
\mathbb{R} \ \ \t{for all} \ \ n\geq 1 \ \ (\t{see eq.}~(13)),
\end{align*}
with $\xi_{5}>0$ independent of $n\geq 1$. So we may assume that
\begin{align*}
u_{n}\xrightarrow{w^{*}} u \quad \t{in} \ \ L^{\infty}(\mathbb{R},\mathbb{R}^{N})
\ \ \t{and} \ \ u_{n}&\xrightarrow{w} u \quad \t{in} \ \ L^{q}(T_{m},\mathbb{R}^{N})\\
&\quad
\t{for all} \ \ m\geq 1 \ \left(\frac{1}{p}+\frac{1}{q}=1\right).
\end{align*}

Evidently $u\in L^{\infty}(\mathbb{R},\mathbb{R}^{N})\cap
L^{q}_{\t{loc}}(\mathbb{R},\mathbb{R}^{N})$ and using Proposition
VII.3.13, p.~694, of [14], we have
$u(t)\in\partial j(t,x(t))$ a.e. on $T_{n}$ for all $n\geq 1$,
hence $u(t)\in\partial j(t,x(t))$ a.e. on $\mathbb{R}$ (recall
that the multifunction $x\rightarrow \partial j(t,x(t))$ is upper
semicontinuous, see [5], p.~29). For every $\tau>0$ we have
that
\begin{align*}
\ &\int_{-\tau}^{\tau}\|x_{n}(t)-x(t)\|^{p}\d t\rightarrow 0 \quad \t{as} \ \ n\rightarrow\infty,\\
\Rightarrow &\mathop {\lim }\limits_{n \to \infty }
\int_{-\tau}^{\tau}\|x_{n}(t)\|^{p}\d t=\int_{-\tau}^{\tau}\|x(t)\|^{p}\d t.
\end{align*}

We can find $n_{0}\geq 1$ such that for all $n\geq n_{0}$ we have
$[-\tau,\tau]\subseteq T_{n_{0}}$ and then using (13) we have
\begin{align*}
\ & \int_{-\tau}^{\tau}\|x_{n}(t)\|^{p}\d t \leq \int_{-n_{0}b}^{n_{0}b}\|x_{n}(t)\|^{p}\d t
\leq\xi_{3}^{p},\\[.5pc]
\Rightarrow & \int_{-\tau}^{\tau}\|x(t)\|^{p}\d t \leq\xi_{3}^{p}.
\end{align*}

Because $\tau >0$ was arbitrary it follows that $x\in L^{p}(\mathbb{R},\mathbb{R}^{N})$.

Next let $\theta\in C_{0}^{\infty}(\mathbb{R},\mathbb{R}^{N})$.
Then supp $\theta \subseteq T_{n}=[-nb,nb]$ for some $n\geq 1$.
Integrating by parts we have $|\int_{\mathbb{R}}
(x_{n}(t),\theta'(t))_{\mathbb{R}^{N}}\d t|=|\int_{\mathbb{R}}
(x_{n}'(t),\theta(t))_{\mathbb{R}^{N}}\d t|$, hence we have
\begin{equation*}
\hskip -2.5pc \left|\int_{\mathbb{R}}
(x_{n}(t),\theta'(t))_{\mathbb{R}^{N}}\d t\right|=\left|\int_{-nb}^{nb}
(x_{n}'(t),\theta(t))_{\mathbb{R}^{N}}\d t\right|\leq\|x'_{n}\|_{L_{n}^{p}}\|\theta\|_{L_{n}^{q}}
\leq \xi_{3}\|\theta\|_{L^{q}(\mathbb{R},\mathbb{R}^{N})}
\end{equation*}

$\left.\right.$\vspace{-1.5pc}

\noindent (see (12) and recall that $\theta\in
C_{0}^{\infty}(\mathbb{R},\mathbb{R}^{N})$). Note that
$(x_{n}(t),\theta'(t))_{\mathbb{R}^{N}}\rightarrow
(x(t),\theta'(t))_{\mathbb{R}^{N}}$ uniformly on compact sets
(i.e. the convergence is in
$C_{\t{loc}}(\mathbb{R},\mathbb{R}^{N})$) and
$|(x_{n}(t),\theta'(t))_{\mathbb{R}^{N}}|\leq
\|x_{n}\|_{L_{n}^{\infty}}\|\theta'(t)\|\leq
\xi_{4}\|\theta'(t)\|$ a.e. on $T_{n}$ (see (13)).

Set
\begin{equation*}
\eta(t)= \begin{cases}
  \xi_{4}\|\theta'(t)\|,& \mbox{if $t\in \t{supp}\ \theta$} \\[2pt]
   0,& \mbox{otherwise}
\end{cases}.
\end{equation*}

Then $\eta\in L^{1}(\mathbb{R})$ and we have
$|(x_{n}(t),\theta'(t))_{\mathbb{R}^{N}}|\leq \eta(t)$ a.e. on
$\mathbb{R}$. By the dominated convergence theorem we have
\begin{align*}
\ & \int_{\mathbb{R}}
(x_{n}(t),\theta'(t))_{\mathbb{R}^{N}}\d t\rightarrow
\int_{\mathbb{R}}
(x(t),\theta'(t))_{\mathbb{R}^{N}}\d t\\[.5pc]
\Rightarrow & \left|\int_{\mathbb{R}}
(x(t),\theta'(t))_{\mathbb{R}^{N}}\d t\right|\leq
\xi_{3}\|\theta\|_{L^{q}(\mathbb{R},\mathbb{R}^{N})}.
\end{align*}

From Proposition IX.3, p.~153 of [3], we obtain that $x\in
W^{1,p}(\mathbb{R},\mathbb{R}^{N})$. Also since
$u_{n}\xrightarrow{w} u$ in $L^{q}_{\t{loc}}
(\mathbb{R},\mathbb{R}^{N})$, we have that
$\int_{\mathbb{R}}(u_{n}(t),\theta(t))_{\mathbb{R}^{N}}
\rightarrow \int_{\mathbb{R}}(u(t),\theta(t))_{\mathbb{R}^{N}}$,
while from the fact that $x_{n}\rightarrow x$ in
$C_{\t{loc}}(\mathbb{R},\mathbb{R}^{N})$ it follows that
\begin{equation*}
\int_{\mathbb{R}}g(t)\|x_{n}(t)\|^{p-2}(x_{n}(t),\theta(t))_{\mathbb{R}^{N}}\d t\rightarrow
\int_{\mathbb{R}}g(t)\|x(t)\|^{p-2}
(x(t),\theta(t))_{\mathbb{R}^{N}}\d t.
\end{equation*}

Also from integration by parts we have
\begin{equation*}
\int_{\mathbb{R}}(\|x'_{n}(t)\|^{p-2}x'_{n}(t),\theta'(t))_{\mathbb{R}^{N}}\d t=
-\int_{\mathbb{R}}((\|x'_{n}(t)\|^{p-2}x_{n}'(t))',\theta(t))_{\mathbb{R}^{N}}\d t.
\end{equation*}

Because $x_{n}$ is a solution of (7), we see that
$(\|x_{n}'(\cdot)\|^{p-2}x_{n}' (\cdot))'\in L^{q}_{n}$ for
all $n\geq 1$. From this we obtain that
$\|x_{n}'(\cdot)\|^{p-2} x_{n}'(\cdot)\in W^{1,q}_{n}$ for
all $n\geq 1$. Also by passing to a subsequence if necessary, we
may assume that
$\|x_{n}'(\cdot)\|^{p-2}x_{n}'(\cdot)\xrightarrow{w} v$ in
$W^{1,q}_{\t{loc}} (\mathbb{R},\mathbb{R}^{N})$, hence
$\|x_{n}'(\cdot)\|^{p-2}x_{n}'(\cdot)\rightarrow v$ in $L^{1}
_{\t{loc}}(\mathbb{R},\mathbb{R}^{N})$ (and in
$C_{\t{loc}}(\mathbb{R},\mathbb{R}^{N})$ too). If $\sigma\!\!:
\mathbb{R}^{N}\rightarrow \mathbb{R}^{N}$ is defined by
$\sigma(x)=\|x\|^{p-2}x$, we have
$\sigma^{-1}(\|x_{n}'(\cdot)\|^{p-2}x_{n}'(\cdot))=x_{n}'(\cdot)\rightarrow
\sigma^{-1}(v)$ in $L_{\rm loc}^{1}(\mathbb{R},\mathbb{R}^{N})$ and
so $\sigma^{-1}(v)=x$, hence
$v(\cdot)=\|x'(\cdot)\|^{p-2}x'(\cdot)$. Therefore we have
\begin{align*}
&\int_{\mathbb{R}}(\|x'_{n}(t)\|^{p-2}x'_{n}(t),\theta'(t))_{\mathbb{R}^{N}}\d t
 \rightarrow \int_{\mathbb{R}}(\|x'(t)\|^{p-2}x'(t),\theta'(t))_{\mathbb{R}^{N}}\d t\\[.3pc]
\ &\quad =-\int_{\mathbb{R}}((\|x'(t)\|^{p-2}x'(t))',\theta'(t))_{\mathbb{R}^{N}}\d t
\ \ (\t{by integration by parts}).
\end{align*}

Since for all $n\geq 1$ large we have
\begin{align*}
&\int_{\mathbb{R}}(\|x'_{n}(t)\|^{p-2}x'_{n}(t),\theta'(t))_{\mathbb{R}^{N}}\d t\\[.3pc]
&\quad +
\int_{\mathbb{R}}g(t)\|x_{n}(t)\|^{p-2}(x_{n}(t),\theta(t))_{\mathbb{R}^{N}}\d t
=\int_{\mathbb{R}}(u(t),\theta(t))_{\mathbb{R}^{N}} \d t
\end{align*}
by passing to the limit as $n\rightarrow \infty$ and using the convergences established above,
we obtain
\begin{align*}
&-\int_{\mathbb{R}}((\|x'(t)\|^{p-2}x'(t))',\theta(t))_{\mathbb{R}^{N}}\d t\\[.3pc]
&\quad +
\int_{\mathbb{R}}g(t)\|x(t)\|^{p-2}(x(t),\theta(t))_{\mathbb{R}^{N}}\d t
=\int_{\mathbb{R}}(u(t),\theta(t))_{\mathbb{R}^{N}} \d t.
\end{align*}

Because $\theta\in C_{0}^{\infty}(\mathbb{R},\mathbb{R}^{N})$ is arbitrary, it follows that
\begin{equation*}
-(\|x'(t)\|^{p-2}x'(t))'+g(t)\|x(t)\|^{p-2}x(t)=u(t) \quad \t{a.e. on} \ \ \mathbb{R}
\end{equation*}
and $u\in L_{\t{loc}}^{q}(\mathbb{R},\mathbb{R}^{N})$, $u(t)\in
\partial j(t,x(t))$ a.e. on $\mathbb{R}$.

Next we show that $x(\pm\infty)=x'(\pm\infty)=0$. Recall that from previous arguments we have
 $x\in W^{1,p}(\mathbb{R},\mathbb{R}^{N})$. So from Corollary VII.8, p.~130 of [3],
we have $x(t)\rightarrow 0$ as $|t|\rightarrow\infty$. Hence we have $x(\pm\infty)=0$.

Since $u(t)\in \partial j(t,x(t))$ a.e. on $\mathbb{R}$, from
hypothesis H(j)$_{3}$(iii) we have $\|u(t)\|\leq
a_{1}(t)(1+\|x(t)\|^{p-1})$ a.e. on $\mathbb{R}$. Because $x\in
W^{1,p}(\mathbb{R},\mathbb{R}^{N})$, we have
$\|x(\cdot)\|^{p-2}x(\cdot)\in L^{q} (\mathbb{R},\mathbb{R}^{N})$
and so $u\in L^{q}(\mathbb{R},\mathbb{R}^{N})$. Therefore,
$\|x'(\cdot)\|^{p-2}x'(\cdot)\in
W^{1,q}(\mathbb{R},\mathbb{R}^{N})$ and once again from
p.~130 of [3], we have that $\|x'(t)\|^{p-1}\rightarrow 0$ as
$|t|\rightarrow\infty$, hence $x'(t)\rightarrow 0$ as
$|t|\rightarrow\infty$. Therefore $x'(\pm\infty)=0$ and we have
proved that $x$ is a homoclinic (to 0) solution.

It remains to show that $x$ is non-trivial. For every $n\geq 1$, we have
\begin{align*}
\ & \;A(x_{n})+g\|x_{n}\|^{p-2}x_{n}=u_{n}\\
\Rightarrow &\; c \|x_{n}\|^{p}_{L^{p}_{n}}\leq\int_{-nb}^{nb}
(u_{n}(t),x_{n}(t))_{\mathbb{R}^{N}}\d t.
\end{align*}

Set
\begin{equation*}
h_{n}(t)=
\begin{cases}
  \frac{(u_{n}(t),x_{n}(t))_{\mathbb{R}^{N}}}{\|x_{n}(t)\|^{p}},& \mbox{if $x_{n}(t)\neq 0$} \\[3pt]
   0,& \mbox{if $x_{n}(t)=0$}
\end{cases}.
\end{equation*}
We have
\begin{equation*}
\hskip -3pc c \|x_{n}\|^{p}_{L^{p}_{n}}\leq\int_{-nb}^{nb}
(u_{n}(t),x_{n}(t))_{\mathbb{R}^{N}}\d t
=\int_{-nb}^{nb}h_{n}(t)\|x_{n}(t)\|^{p} \d t\leq \mathop
{\t{ess}\sup }\limits_{T_n} h_{n}\|x_{n}\|_{L^{p}_{n}}^{p}.
\end{equation*}

By virtue of hypothesis $\hbox{H(j)}_{3}\hbox{(v)}$ (see the remark following
H(j)$_{3}$), given $\varepsilon >0$ we can find $\delta=\delta
(\varepsilon)>0$ such that for almost all $t\in \mathbb{R}$, all
$\|x\|\leq \delta$ and all $u\in \partial j(t,x)$ we have
\begin{equation}
\frac{(u,x)_{\mathbb{R}^{N}}}{\|x\|^{p}}\leq \varepsilon.
\end{equation}

If $x=0$, then $x_{n}\rightarrow 0$ in
$C_{\t{loc}}(\mathbb{R},\mathbb{R}^{N})$ and so we can find
$n_{0}\geq 1$ such that for all $n\geq n_{0}$ and all $t\in
T_{n}$, we have $\|x_{n}(t)\|\leq \delta$. Therefore for all
$n\geq n_{0}$ and almost all $t\in T_{n}$, we have $h_{n}(t)\leq
\varepsilon$ (see (14)) and so $c\leq \mathop {\t{ess}\sup
}_{T_{n}}h_{n}= \mathop {\t{ess}\sup }_\mathbb{R}
h_{n}\leq \varepsilon$ for all $n\geq n_{0}$ (recall that
$x_{n},u_{n}$ were extended by periodicity to all of
$\mathbb{R}$). Let $\varepsilon\downarrow 0$ to obtain $0<c \leq
0$, a contradiction. This proves that $x\neq 0$.

Therefore $x\in C(\mathbb{R},\mathbb{R}^{N})\cap W^{1,p}(\mathbb{R},\mathbb{R}^{N})$ is the
desired non-trivial, homoclinic (to 0) solution of the non-smooth non-linear periodic system.\hfill QED
\end{proof}

\section{Scalar equations}

In this last part of the paper we study the scalar (i.e. $N=1$) problem. We approach the problem
using a generalized Landesmann--Lazer type condition, which is more general than the one used by Tang
[33] in the context of smooth semilinear periodic equations. So our work is a two-fold
generalization of the work of Tang.

First we examine the following non-linear scalar periodic problem:
\begin{equation}
\left \{  \begin{array}{l} -(|x'(t)|^{p-2}x'(t))'
\in \partial j(t,x(t)) \quad  \t{a.e} \ \t{on} \ T=[0,b]\\
x(0)=x(b),\quad  x'(0)=x'(b),\quad 1< p<\infty
\end{array} \right \}.
\end{equation}

The conditions on the non-smooth potential $j$ are the following:
\begin{itemize}
\leftskip 1pc
\item[{H(j)$_4$:}] $j\!\!:\,T\times\mathbb{R}\longmapsto \mathbb{R}$ is a
functional such that $j(\cdot,0)\in L^{1}(T)$ and

%\begin{itemize}
\hskip 1pc
 \leftskip .8pc
\item[(i)] for all $x\in \mathbb{R}$, $t\longmapsto j(t,x)$ is
measurable;

\item[(ii)] for almost all $ t\in T$, the function $x\longmapsto j(t,x)$
is locally Lipschitz;

\item[(iii)] for almost all $t\in T$, all $x\in \mathbb{R}$ and all
$u\in \partial j(t,x)$, we have
\begin{equation*}
\hskip -3.2pc|u|\le a_{1}(t)+c_{1}(t)|x|^{r-1},
\end{equation*}
$1\le r<+\infty$ with $a_{1},c_{1}\in L^{r'}(T)$,
$\frac{1}{r}+\frac{1}{r'}=1$;

\item[(iv)] $\mathop {\lim }\limits_{|x| \to \infty }\frac{u}{x}=0$
uniformly for almost all $t\in T$ and all $u\in\partial j(t,x)$;

\item[(v)] there exist functions $j_{+},j_{-}\in L^{1}(T)$ such that
$j_{+}(t)=\mathop {\lim \inf }\limits_{x \to +\infty }\frac{j(t,x)}{x}$
and $j_{-}(t)=\mathop {\lim \sup}\limits_{x \to -\infty
}\frac{j(t,x)}{x}$ a.e. on $T$ and $\int_{0}^{b}j_{-}(t)\d
t<0<\int_{0}^{b}j_{+}(t)\d t$.
%\end{itemize}
\end{itemize}

We consider the locally Lipschitz functional $\varphi\!\!:W_{\rm per}^{1,p}(T)\rightarrow\mathbb{R}$
defined by
\begin{equation*}
\varphi(x)=\frac{1}{p}\|x'\|_{p}^{p}-\int_{0}^{b}j(t,x(t))\d t.
\end{equation*}

\setcounter{prop}{3}

\begin{prop}$\left.\right.$\vspace{.5pc}

\noindent If hypotheses $\hbox{\rm H(j)}_{4}$ hold{\rm ,}
then $\varphi$ satisfies the non-smooth PS-condition.
\end{prop}

\begin{proof}
We consider a sequence $\{x_{n}\}_{n\geq 1}\subseteq W_{\rm per}^{1,p}(T)$
such that
\begin{equation*}
|\varphi(x_{n})|\leq M_{1} \quad \t{for some} \ \ M_{1}>0 \ \ \t{and all} \ \ n\geq 1 \ \
\t{and} \ \ m(x_{n})\rightarrow 0.
\end{equation*}

As before we choose $x_{n}^{*}\in \partial \varphi(x_{n})$ such that $m(x_{n})=\|x_{n}^{*}\|$,
$n\geq 1$. We have
\begin{equation*}
x_{n}^{*}=A(x_{n})-u_{n} \quad \t{with} \ \ u_{n}\in L^{r'}(T),u_{n}(t)\in\partial j(t,x_{n}(t))
\ \ \t{a.e. on} \ \ T.
\end{equation*}

We claim that $\{x_{n}\}_{n\geq 1}\subseteq W_{\rm per}^{1,p}(T)$ is bounded. Suppose that this is
not the case. By passing to a subsequence if necessary, we may assume that
$\|x_{n}\|\rightarrow\infty$. Let $y_{n}=\frac{x_{n}}{\|x_{n}\|}$, $n\geq 1$. We may assume that
\begin{equation*}
y_{n}\xrightarrow{w} y \quad \t{in} \ \ W_{\rm per}^{1,p}(T) \quad \t{and} \quad y_{n}\rightarrow y \ \
\t{in} \ \ C_{\rm per}(T).
\end{equation*}
(Recall that $W_{\rm per}^{1,p}(T)$ is embedded compactly in $C_{\rm per}(T)$.) From the choice of the
sequence $\{x_{n}\}_{n\geq 1}\subseteq W_{\rm per}^{1,p}(T)$, we have
\begin{equation}
\frac{|\varphi(x_{n})|}{\|x_{n}\|^{p}}=\left|\frac{1}{p}\|y_{n}'\|_{p}^{p}-\int_{0}^{b}
\frac{j(t,x_{n}(t))}{\|x_{n}\|^{p}}\d t\right|\leq\frac{M_{1}}{\|x_{n}\|^{p}}.
\end{equation}

By virtue of hypothesis H(j)$_{4}$(v) we have $\int_{0}^{b}
\frac{j(t,x_{n}(t))}{\|x_{n}\|^{p}}\d t\rightarrow 0$. So from (16)
and the weak lower semicontinuity of the norm in a Banach space,
we obtain $\frac{1}{p}\|y'\|_{p}^{p}=0$, hence
$y=c\in\mathbb{R}$. If $c=0$, then we have $\|y'\|_{p}\rightarrow
0$, hence $y_{n}\rightarrow 0$ in $W_{\rm per}^{1,p}(T)$. But for
every $n\geq 1$, $\|y_{n}\|=1$ and so we have a contradiction.
Therefore $y=c\neq 0$ and without any loss of generality we may
assume that $y=c>0$ (the analysis is the same if instead we
assume that $y=c<0$). Recall that
$W_{\rm per}^{1,p}(T)=\mathbb{R}\oplus V$ with $V=\{v\in
W^{1,p}_{\rm per}(T)\!\!:\int_0^bv(t)\d t=0\}$. We have
$x_{n}=\bar{x}_{n}+\hat{x}_{n}$ with $\bar{x}_{n}\in\mathbb{R}$,
$\hat{x}_{n}\in V$, $n\geq 1$. Then
$y_{n}=\bar{y}_{n}+\hat{y}_{n}$ with
$\bar{y}_{n}=\frac{\bar{x}_{n}}{\|x_{n}\|}$,
$\hat{y}_{n}=\frac{\hat{x}_{n}}{\|x_{n}\|}$, $n\geq 1$. From the
choice of the sequence $\{x_{n}\}_{n\geq 1}\subseteq
W_{\rm per}^{1,p}(T)$, we have $-\varepsilon_{n}\leq \langle x_{n}^{*},
y_{n}\rangle \leq\varepsilon_{n}$ with $\varepsilon_{n}\downarrow 0$,
hence
\begin{equation}
-\varepsilon_{n}\leq \frac{1}{\|x_{n}\|}\left[\|\hat{x}_{n}'\|_{p}^{p}-\int_{0}^{b}
u_{n}(t)\hat{x}_{n}(t)\d t\right]\leq\varepsilon_{n}.
\end{equation}

Since we have assumed that $y=c>0$, we have that for all $t\in T$,
$x_{n}(t)\rightarrow +\infty$ as $n\rightarrow\infty$. We claim
that this convergence is uniform in $t\in T$. Indeed let
$\varepsilon>0$ be such that $0<\varepsilon<c$ (recall that we
have assumed $c>0$). Since $y_n\to c$ in $C(T)$, we can find
$n_0\geq 1$ such that for all $n\geq n_0$ and all $t\in
T,\;|y_n(t)-c|<\varepsilon$, hence $0<c-\varepsilon<|y_n(t)|$.
Because $\|x_n\|\to\infty$ given $\beta>0$ we can find $n_1\geq
1$ such that for all $n\geq n_1$ we have $\|x_n\|\geq\beta>0$. So
for all $n\geq n_2=\max\{n_0,n_1\}$ and all $t\in T$ we have
$\frac{|x_n(t)|}{\beta}\geq\frac{|x_n(t)|}{\|x_n\|}=|y_n(t)|>c-\varepsilon=\theta
>0\Rightarrow |x_n(t)|\geq\beta\theta>0$. Because $\beta>0$ is
arbitrary and $\theta>0$, we infer that $\min_{T}
|x_n(t)|\to +\infty$. Then we have
\begin{equation*}
\int_{0}^{b}u_{n}(t)\hat{x}_{n}(t)\d t=\int\limits_{\{ x_n (t) \ne 0\} } {\frac{u_{n}(t)}{x_{n}(t)}
x_{n}(t)\hat{x}_{n}(t)}\d t.
\end{equation*}

Evidently $|\{x_{n}= 0\}|\rightarrow 0$, while by virtue of
hypothesis H(j)$_{4}$(iv), given $\varepsilon >0$ we can find
$n_{3}\geq 1$ such that for all $n\geq n_{3}$ and almost all $t\in
T$ we have $|{u_{n}(t)}/{x_{n}(t)}|<\varepsilon$. So for
$n\geq n_{3}$ we have
\begin{equation*}
\left|\,\int\limits_{\{ x_n (t) \ne 0\} } {\frac{u_{n}(t)}{x_{n}(t)}
x_{n}(t)\hat{x}_{n}(t)}\d t\right|\leq\varepsilon \|\hat{x}_{n}\|_{2}^{2}\leq\varepsilon c_{2}
\|\hat{x}_{n}\|^{p} \quad \t{for some} \ \ c_{2}>0.
\end{equation*}

From (17) and the Poincar\'e--Wirtinger inequality, we have
\begin{equation*}
 \frac{1}{\|x_{n}\|}[c_{3}\|\hat{x}_{n}\|^{p}-\varepsilon c_{2}
\|\hat{x}_{n}\|^{p}]=(c_{3}-\varepsilon c_{2})
\frac{\|\hat{x}_{n}\|^{p}}{\|x_{n}\|}\leq\varepsilon_{n} \quad \t{for some} \ \ c_{3}>0.
\end{equation*}

Choose $\varepsilon >0$ such that $\varepsilon c_{2}<c_{3}$. Since $\varepsilon_{n}\downarrow 0$
it follows that $\frac{\|\hat{x}_{n}\|^{p}}{\|x_{n}\|}\rightarrow 0$, hence once again by the
Poincar\'e--Wirtinger inequality, we have
\begin{equation}
\frac{\|x'_{n}\|^{p}_{p}}{\|x_{n}\|}\rightarrow 0 \quad \t{as} \ \ n\rightarrow\infty.
\end{equation}

Recall that for $n\geq 1$, we have
\begin{equation}
-\frac{M_{1}}{\|x_{n}\|}\leq \frac{1}{p}\frac{\|x'_{n}\|_{p}^{p}}{\|x_{n}\|}-
\int_{0}^{b}\frac{j(t,x_{n}(t))}{\|x_{n}\|}\d t\leq\frac{M_{1}}{\|x_{n}\|}.
\end{equation}
(We assume without any loss of generality that $\|x_{n}\|\geq\xi >0$ for all $n\geq 1$; recall
that $\|x_{n}\|\rightarrow\infty$.) We can write that
\begin{equation*}
\int_{0}^{b}\frac{j(t,x_{n}(t))}{\|x_{n}\|}\d t+
\int\limits_{\{ x_n (t) \ne 0\} }\frac{j(t,x_{n}(t))}{x_{n}(t)}y_{n}(t)\d t+
\int\limits_{\{ x_n (t) = 0\} }\frac{j(t,0)}{\|x_{n}\|}\d t.
\end{equation*}

By virtue of hypothesis H(j)$_{4}$(v), given $\varepsilon >0$ we
can find $M=M(\varepsilon)>0$ such that for almost all $t\in T$
and all $x\geq M$, we have ${j(t,x)}/{x}\geq j_{+}(t)
-\varepsilon$. Recall that $x_{n}(t)\rightarrow +\infty$ uniformly
in $t\in T$ (i.e. $\mathop {\min}_T x_{n}\rightarrow
+\infty)$. So we can find $n_{0}\geq 1$ such that for all $n\geq
n_{0}$ we have
\begin{align*}
\ & j_{+}(t)-\varepsilon \leq\frac{j(t,x_{n}(t))}{x_{n}(t)} \quad \t{a.e. on} \ \ T\\
\Rightarrow &\int\limits_{\{ x_n (t) \ne 0\} }(j_{+}(t)-\varepsilon )y_{n}(t)\leq
\int\limits_{\{ x_n (t) \ne 0\} }\frac{j(t,x_{n}(t))}{x_{n}(t)}y_{n}(t)\d t\\
\Rightarrow & \int_{0}^{b}j_{+}(t)c \d t\leq\mathop {\lim \inf }\limits_{n \to \infty }
\int\limits_{\{ x_n (t) \ne 0\} }\frac{j(t,x_{n}(t))}{x_{n}(t)}y_{n}(t)\d t \\
&\qquad\qquad\qquad\qquad\qquad (\t{since
$\varepsilon >0$ was arbitrary}).
\end{align*}

Also we have $\int_{\{ x_n (t)= 0\}} {(j(t,0)}/{\|x_{n}\|})\d t\rightarrow 0$. So finally we have
\begin{equation*}
c\int_{0}^{b}j_{+}(t) \d t\leq\mathop {\lim \inf }\limits_{n \to \infty }
\int_{0}^{b}\frac{j(t,x_{n}(t))}{\|x_{n}\|}\d t.
\end{equation*}

Using this and (18) in (19) we obtain $c\int_{0}^{b}j_{+}(t)
\d t\leq 0$, hence $\int_{0}^{b}j_{+}(t) \d t\leq 0$ (recall that we
have assumed that $c>0$). This contradicts hypothesis
H(j)$_{4}$(v). Therefore $\{x_{n}\}_{n\geq 1}\subseteq
W_{\rm per}^{1,p}(T)$ is bounded and then arguing as in the proof of
Theorem 1 we conclude that $\varphi$ satisfies the non-smooth
PS-condition.\hfill QED
\end{proof}

\begin{prop}$\left.\right.$\vspace{.5pc}

\noindent If hypotheses $\hbox{\rm H(j)}_{4}$ hold{\rm ,}
then $\mathop {\lim}_{|\mathop {\xi | \to \infty }
\limits_{\xi  \in \mathbb{R}} }\varphi(\xi)=-\infty $ {\rm (}i.e. $\varphi|_{\mathbb{R}}$ is
anticoercive{\rm )}.
\end{prop}

\begin{proof}
Suppose that the result of the proposition is not true. Then we can
find $\{\xi_{n}\}_{n\geq 1}\subseteq\mathbb{R}$ such that $|\xi_{n}|\rightarrow\infty$ and
$\gamma\in \mathbb{R}$ such that $\gamma\leq\varphi(\xi_{n})$ for all $n\geq 1$. So $\gamma\leq
\mathop {\lim \inf}_{n \to \infty }\varphi(\xi_{n})=
\mathop {\lim \inf}_{n \to \infty }(-\int_{0}^{b}j(t,\xi_{n})\d t)$. First suppose that
$\xi_{n}\rightarrow +\infty$. We may assume that $\xi_{n}>0$ for all $n\geq 1$. We have
\begin{align*}
\ & \frac{\gamma}{\xi_{n}}\leq-\int_{0}^{b}\frac{j(t,\xi_{n})}{\xi_{n}}\d t,\\
\Rightarrow &\mathop {\lim \sup }\limits_{n \to \infty }\int_{0}^{b}\frac{j(t,\xi_{n})}
{\xi_{n}}\d t\leq 0.
\end{align*}

On the other hand as in the proof of Proposition~4, we obtain that
\begin{equation*}
\int_{0}^{b}j_{+}(t)\d t \leq\mathop {\lim \inf }\limits_{n \to \infty }\int_{0}^{b}\frac{j(t,\xi_{n})}
{\xi_{n}}\d t\leq 0,
\end{equation*}
which is a contradiction. Similarly if $\xi_{n}\rightarrow -\infty$, we obtain
$\int_{0}^{b}j_{-}(t)\d t\geq 0 $, again a contradiction.\hfill QED
\end{proof}

Recall the direct sum decomposition $W_{\rm per}^{1,p}(T)=\mathbb{R}\oplus V$ with
$V=\{v\in W^{1,p}_{\rm per}(T)\!\!:\int_0^bv(t)\d t=0\}$.

\begin{prop}$\left.\right.$\vspace{.5pc}

\noindent If hypotheses $\hbox{\rm H(j)}_{4}$ hold{\rm ,}
then $\varphi|_{V}$ is coercive {\rm (}i.e.
$\varphi(v)\rightarrow + \infty$ as $\|v\|\rightarrow\infty, v\in
V${\rm )}.
\end{prop}

\begin{proof}
For every $v\in V$ we have $\varphi(v)=\frac{1}{p}\|v'\|_{p}^{p}-
\int_{0}^{b}j(t,v(t))\d t$. From the Poincar\'e--Wirtinger inequality we know that on $V$,
$\|v'\|_{p}$ is an equivalent norm. So
\begin{equation*}
\frac{\varphi(v)}{\|v\|^{p}}\geq c_{4}-\int_{0}^{b}\frac{j(t,v(t))}{\|v\|^{p}}\d t \quad \t{for some}
\ \ c_{4}>0.
\end{equation*}

Recall that $\int_{0}^{b}({j(t,v(t))}/{\|v\|^{p}})\d t \rightarrow 0$. So
$\mathop {\lim \inf}_{\mathop {||v|| \to \infty }\limits_{v \in V} }\frac{\varphi(v)}
{\|v\|^{p}}\geq c_{4}>0$, hence $\varphi(v)\rightarrow +\infty$ as $\|v\|\rightarrow\infty$,
$v\in V$.\hfill QED
\end{proof}

\setcounter{lem}{6}

These auxiliary results lead to the following existence theorem.
\begin{thm}[\!]
If hypotheses $\hbox{\rm H(j)}_{4}$ hold{\rm ,}
then problem $(15)$ has at least one solution
$x\in C^{1}(T)$ with $|x'(\cdot)|^{p-2}x'(\cdot)\in W^{1,r'}(T)$.
\end{thm}

\begin{proof}
Propositions~4--6 permit the
application of the non-smooth saddle point theorem. So we obtain
$x\in W^{1,p}_{\rm per}(T)$ such that $0\in\partial\varphi(x)$, hence
$A(x)=u$ with $u\in L^{r'}(T)$, $u(t)\in \partial j(t,x(t))$ a.e.
on $T$. As in the proof of Theorem 1 we show that $x\in C^{1}(T)$,
$|x'(\cdot)|^{p-2}x'(\cdot)\in W^{1,r'}(T)$ and it solves problem
(15).\hfill QED
\end{proof}

As we have already mentioned in the beginning of this section, our
generalized Landesmann--Lazer type condition (see hypothesis
H(j)$_{4}$(v)) generalizes the one used by Tang [33] (for smooth
potentials). In the next proposition we are going to show this.
For this purpose we introduce the following functions:
\begin{align*}
g_{1}(t,x)&=\min\{u\!\!:u\in\partial j(t,x)\} \quad \t{and} \quad
g_{2}(t,x)=\max\{u\!\!:u\in\partial j(t,x)\},\\
G_{1}(t,x)&=
\begin{cases}
  \frac{2j(t,x)}{x}-g_{1}(t,x),& \mbox{if $|x|\neq 0$} \\
    0,& \mbox{if $|x|=0$}
\end{cases}
\end{align*}
and
\begin{align*}
G_{2}(t,x)&=
\begin{cases}
  \frac{2j(t,x)}{x}-g_{2}(t,x),& \mbox{if $|x|\neq 0$} \\
    0,& \mbox{if $|x|=0$}
\end{cases},\\
G_{1}^{-}(t)&=\mathop {\lim \sup }\limits_{x \to  - \infty }G_{1}(t,x)\quad \hbox{and}
\quad G_{2}^{+}(t)=\mathop {\lim \inf }\limits_{x \to  + \infty }G_{2}(t,x).
\end{align*}

The functions $G_{1}^{-}$, $G_{2}^{+}$ are essentially the ones
used by Tang [33] in the context of smooth, semilinear (i.e.
$p=2$) periodic problems. In that case, since $j(t,\cdot)\in
C^{1}(\mathbb{R})$, we have $g_{1}=g_{2}$, and hence $G_{1}=G_{2}$.

\setcounter{prop}{7}

\begin{prop}$\left.\right.$\vspace{.5pc}

\noindent For all $t\in T\setminus D$ with $|D|=0${\rm ,} $G_{2}^{+}(t)\leq j_{+}(t)$
and $j_{-}(t)\leq G_{1}^{-}(t)$.
\end{prop}

\begin{proof}
Let $D\subseteq T$ the Lebesgue-null
set outside of which hypotheses \hbox{$\hbox{H(j)}_{4}\hbox{(ii)}\rightarrow (v)$}
hold. Let $t\in T\setminus D$ and set $k_{\varepsilon}^{+}(t)=
G_{2}^{+}(t)-\varepsilon$. We can find $M_{1}>0$ such that for
all $x\geq M_{1}$ we have
\begin{align*}
\ & G_{2}^{+}(t)-\varepsilon = k_{\varepsilon}^{+}(t)\leq G_{2}(t,x)\\
\Rightarrow & \frac{k_{\varepsilon}^{+}(t)}{x^{2}}= \frac{\d}{\d x}\left(-\frac{k_{\varepsilon}^{+}(t)}
{x}\right)\leq\frac{G_{2}(t,x)}{x^{2}}.
\end{align*}

For all $u\in\partial j(t,x)$, we have
\begin{equation*}
\frac{G_{2}(t,x)}{x^{2}}=\frac{2j(t,x)}{x^{3}}-\frac{g_{2}(t,x)}{x^{2}}
\leq\frac{2j(t,x)}{x^{3}}-\frac{u}{x^{3}}.
\end{equation*}

From p.~48 of [5], we know that $x\rightarrow {2j(t,x)}/{x^{2}}$ is locally Lipschitz
on $[M_1,+\infty)$ and
\begin{equation*}
\partial\left(\frac{j(t,x)}{x^{2}}\right)\subseteq\frac{\partial j(t,x)}{x^{2}}-\frac{2j(t,x)}{x^{3}}.
\end{equation*}

Therefore for all $t\in T$, all $x\geq M_{1}$ and all $u\in \partial j(t,x)$, we have
\begin{align*}
\ & u\leq\frac{g_{2}(t,x)}{x^{2}}-\frac{2j(t,x)}{x^{3}}=-\frac{1}{x^{2}}G_{2}(t,x)\\
\Rightarrow & u\leq\frac{\d}{\d x}\left(\frac{k_{\varepsilon}^{+}(t)}{x}\right).
\end{align*}

Since for $t\in T\setminus D$, the function $t\rightarrow {j(t,x)}/{x^{2}}$ is locally
Lipschitz on $[M_1,+\infty)$, it is differentiable at every
$x\in[M_1,+\infty)\setminus D_{1}(t)$, $|D_{1}(t)|=0$. We set
\begin{equation*}
u_{0}(t,x)=
\begin{cases}
  \frac{\d}{\d x}(\frac{j(t,x)}{x^{2}}),& \mbox{if $x\in[M_{1},+\infty)\setminus D_{1}(t)$}, \\[2pt]
    0,& \mbox{if $x\in D_{1}(t)$}.
\end{cases}
\end{equation*}

We fix $t\in T\setminus D$ and choose $x\in [M_{2},+\infty)\setminus D_{1}(t)$. Then
$u_{0}(t,x)\in\partial ({j(t,x)}/{x^{p}})$ and so
\begin{equation}
u_{0}(t,x)=\frac{\d}{\d x}\left(\frac{j(t,x)}{x^{2}}\right)\leq \frac{\d}{\d x}\frac{k_{\varepsilon}^{+}(t)}{x}.
\end{equation}

Let $y<x$ and $y\in [M_{1},+\infty)\setminus D_{1}(t)$. We
integrate (20) over the interval $[y,x]$ and obtain
\begin{equation}
\frac{j(t,x)}{x^{2}}-\frac{j(t,y)}{y^{2}}\leq k_{\varepsilon}^{+}(t)\left(\frac{1}{x}-\frac{1}{y}\right).
\end{equation}

By virtue of hypotheses H(j)$_{4}$(iii), (iv) and the Lebourg mean
value theorem, given $\varepsilon >0$ for all $t\in T\setminus D$
and all $x\geq 0$, we have
\begin{align*}
\ &-\varepsilon x^{2}-c_{1}x +j(t,0)\leq j(t,x) \ \ \t{for some} \ \ c_{1}>0\\
\Rightarrow & -\varepsilon\leq\mathop {\lim \inf }\limits_{x \to  + \infty }
\frac{j(t,x)}{x^{2}},\\
\Rightarrow &\;0 \leq\mathop {\lim \inf }\limits_{x \to  + \infty
} \frac{j(t,x)}{x^{2}}.
\end{align*}

So if we go to (21) and pass to the limit as $x\rightarrow
+\infty$, we obtain
\begin{align*}
\ & \;k_{\varepsilon}^{+}(t)\leq\frac{j(t,y)}{y},\\
\Rightarrow & \;G_{2}^{+}(t)\leq\mathop {\lim \inf }\limits_{x
\to  + \infty } \frac{j(t,y)}{y}=j_{+}(t).
\end{align*}

Similarly we obtain that for all $t\in T\setminus D$, $|D|=0$, we
have $j_{-}(t)\leq G_{1}^{-}(t)$.\hfill QED
\end{proof}

\begin{rem}
{\rm This proposition shows that our
generalized Landesman--Lazer type condition (hypothesis
H(j)$_{4}$(v)) is more general than the one used by Tang [33].
Here is an example of a non-smooth locally Lipschitz potential
which satisfies H(j)$_{4}$(v) but does not satisfy the condition
of Tang. Again for simplicity we drop the $t$-dependence
\begin{equation*}
j(x)= \max \{x^{{1}/{3}},x^{{1}/{2}}\}+\ln (1+|x|)+\cos x +x.
\end{equation*}

A simple calculation shows that $j_{+}=1$, $j_{-}=-1$ but $G_{1}^{-}=G_{2}^{+}=0$.}
\end{rem}

When dealing with the semilinear (i.e $p=2$) case, we can consider
problems at resonance in an eigenvalue of any order. Similar
problems (but with smooth potential) were studied by Mawhin and Ward
[25], p.~67 of Mawhin and Willem [26], Mawhin and Schmitt [24] (problems
near resonance) and Tang [33] (who employed his more restrictive
version of the generalized Landesman--Lazer condition (see
Proposition~8).

The problem under consideration is the following:
\begin{equation}
\left \{  \begin{array}{l} -x''(t)-m^{2}\omega^{2}x(t)\
\in \partial j(t,x(t))-h(t) \quad  \t{a.e} \ \t{on} \ T=[0,b]\\
x(0)=x(b)  ,  x'(0)=x'(b),\; h\in L^{1}(T)
\end{array} \right \}.
\end{equation}

Here $m\in \mathbb{N}_{0}=\{0,1,2,...\}$ and $\omega ={2\pi}/{b}$ (see \S2). Our
hypotheses on $j(t,x)$ are the following:

\begin{itemize}
\leftskip 1pc
\item[{H(j)$_5$:}] $j\!\!:\,T\times\mathbb{R}\longmapsto \mathbb{R}$ is
a functional such that $j(\cdot,0)\in L^{1}(T)$ and

%\begin{itemize}
\hskip 1pc
\leftskip .9pc
\item[(i)] for all $x\in \mathbb{R}$, $t\longmapsto j(t,x)$ is
measurable;

\item[(ii)] for almost all $ t\in T$, the function $x\longmapsto j(t,x)$
is locally Lipschitz;

\item[(iii)] for almost all $t\in T$, all $x\in \mathbb{R}$ and all
$u\in \partial j(t,x)$, we have
\begin{equation*}
\hskip -3.2pc|u|\le a_{1}(t)(1+|x|^{r-1}),
\end{equation*}
$1\le r<+\infty$ with $a_{1}\in L^{r'}(T)$,
$\frac{1}{r}+\frac{1}{r'}=1$;

\item[(iv)] $\mathop {\lim }\limits_{|x| \to \infty }\frac{u}{x}=0$
uniformly for almost all $t\in T$ and all $u\in\partial j(t,x)$;

\item[(v)] there exist functions $j_{\pm}\in L^{1}(T)$ such that
$j_{+}(t)=\mathop {\lim \inf }_{x \to +\infty }\frac{j(t,x)}{x}$ and
$j_{-}(t)=\mathop {\lim \sup}_{x \to -\infty }\frac{j(t,x)}{x}$ a.e. on
$T$ and $\int_{0}^{b}h(t)\sin(m\omega t+\theta)\d t
<\int_{0}^{b}j_{+}(t)\sin(m\omega t+\theta)^{+}-j_{-}(t) \sin(m\omega
t+\theta)^{-}\d t$ for all $\theta\in\mathbb{R}$.
%\end{itemize}
\end{itemize}

In our analysis of problem (22) we shall use the following
subspaces of $W_{\rm per}^{1,2}(T)$:
\begin{align*}
\ \bar{H}=&\;\t{span} \{ \sin k \omega t,\cos k \omega t\!: k=0,1,...,m-1\},\\
\ N_{m} =&\;\t{span}\{\sin m\omega t,\cos m\omega t\},\\
\ \hat{H} =&\;(\bar{H}+N_{m})^{\bot}=\t{span}\{\sin k\omega t,\cos
k\omega t\!\!:k\geq m+1\}.
\end{align*}

We have $W_{\rm per}^{1,2}(T)=\bar{H}\oplus N_{m}\oplus \hat{H}$ and so if $x\in W_{\rm per}^{1,2}(T)$,
we have $x=\bar{x}+x^{0}+\hat{x}$ with $\bar{x}\in \bar{H}$, $x^{0}\in N_{m}$ and
$\hat{x}\in \hat{H}$.

We start with an auxiliary result concerning the subspace $\hat{H}$.
\setcounter{lem}{8}
\begin{lem}
There exists $c>0$ such that for all $x\in \hat{H}$ we have $c\|x\|^{2}\leq\|x'\|_{2}^{2}-
\lambda_{m}\|x\|^{2}_{2}$.
\end{lem}

\begin{proof}
Let $\psi(x)=\|x'\|^{2}_{2}-\lambda_{m}\|x\|^{2}_{2}$ and suppose
that the result is not true. We can find $\{x_{n}\}_{n\geq 1}\subseteq\hat{H}$ such that
$\psi(x_{n})\downarrow 0$. Set $y_{n}={x_{n}}/{\|x_{n}\|}$, $n\geq 1$. We may assume that
$y_{n}\xrightarrow{w} y$ in $W_{\rm per}^{1,2}(T)$ and $y_{n}\rightarrow y$ in $L^{2}(T)$. Thus
in the limit as $n\rightarrow\infty$, we obtain $\|y'\|_{2}^{2}\leq \lambda_{m}\|y\|_{2}^{2}$
and so $y=0$. Hence $\|y_{n}'\|_{2}\rightarrow 0$ and so $y_{n}\rightarrow 0$ in
$W_{\rm per}^{1,2}(T)$, a contradiction to the fact that $\|y_{n}\|=1$ for all $n\geq 1$.
\hfill QED
\end{proof}

We introduce the locally Lipschitz functional $\varphi\!\!:W_{\rm per}^{1,2}(T)\rightarrow \mathbb{R}$
defined by
\begin{equation*}
\varphi(x)=\frac{1}{2}\|x'\|_{2}^{2}-\frac{m^{2}\omega^{2}}{2}\|x\|_{2}^{2}-
\int_{0}^{b}j(t,x(t))\d t+\int_{0}^{b}h(t)x(t)\d t.
\end{equation*}

\setcounter{prop}{9}
\begin{prop}$\left.\right.$\vspace{.5pc}

\noindent If hypotheses $\hbox{\rm H}(\hbox{\rm j})_{5}$ hold{\rm ,}
then $\varphi$ satisfies the non-smooth PS-condition.
\end{prop}

\begin{proof}
Let $\{x_{n}\}_{n\geq 1}\subseteq W_{\rm per}^{1,2}(T)$ be a sequence
such that
\begin{equation*}
|\varphi(x_{n})|\leq M_{1} \ \ \t{for all} \ \ n\geq 1 \ \ \t{and some} \ \ M_{1}>0
\end{equation*}
and
\begin{equation*}
m(x_{n})\rightarrow 0 \ \ \t{as} \ \ n\rightarrow\infty.
\end{equation*}

Again we find $x_{n}^{*}\in \partial\varphi(x_{n})$ such that $m(x_{n})=\|x_{n}^{*}\|$ and
$x_{n}^{*}=A(x_{n})-m^{2}\omega^{2}x_{n}-u_{n}+h$, with $u_{n}\in L^{r'}\!(T)$, $u_{n}(t)\in
\partial j(t,x_{n}(t))$ a.e. on $T$. Here $A\in \mathcal{L}(W_{\rm per}^{1,2}(T),W_{\rm per}^{1,2}(T)^{*})$
is defined $\dual{A(x),y}=\int_{0}^{b}x'(t)y'(t)\d t$ for all
$x,y\in W_{\rm per}^{1,2}(T)$. Of course $A$ is a maximal monotone,
bounded linear operator. We claim that $\{x_{n}\}_{n\geq
1}\subseteq W_{\rm per}^{1,2}(T)$ is bounded. If this is not the
case, we may assume that $\|x_{n}\|\rightarrow\infty$. We set
$y_{n}=\frac{x_{n}}{\|x_{n}\|}$, $n\geq 1$ and we may assume that
\begin{equation*}
y_{n}\xrightarrow{w} y\ \ \t{in} \ \ W_{\rm per}^{1,2}(T) \ \ \t{and} \ \ y_{n}\rightarrow y \ \
\t{in} \ \ C_{\rm per}(T) \ \ \t{as} \ \ n\rightarrow\infty.
\end{equation*}

From the choice of the sequence $\{x_{n}\}_{n\geq 1}\subseteq
W_{\rm per}^{1,2}(T)$ we have
\begin{align}
\ & \;|\langle x_{n}^{*},v\rangle |\leq\varepsilon_{n}\|v\| \ \ \t{for all}
\ \ v\in W_{\rm per}^{1,2}(T) \ \
\t{with} \ \ \varepsilon_{n}\downarrow 0,\nonumber\\
\Rightarrow &\; \left|\int_{0}^{b}y_{n}' v'
\d t-\lambda_{m}\int_{0}^{b}y_{n}v \d t-\int_{0}^{b}
\frac{u_{n}}{\|x_{n}\|}v\ \d t +\int_{0}^{b}\frac{h}{\|x_{n}\|}v\
\d t\right|\leq \varepsilon_{n} \frac{\|v\|}{\|x_{n}\|}\nonumber\\
&\hskip 6cm \t{with} \ \
\lambda_{m}=m^{2}\omega^{2}.
\end{align}

By virtue of hypothesis H(j)$_5$(iv), given $\varepsilon >0$ we
can find $M>0$ such that for almost all $t\in T$, all $|x|\geq M$
and all $u\in \partial j(t,x)$, we have
$|{u}/{x}|\leq\varepsilon$. So we can write that
\begin{align*}
\ & \left|\int_{0}^{b}\frac{u_{n}}{\|x_{n}\|}v
\d t\right|=\left|\int_{\{|x_{n}(t)|\geq M\}}\frac{u_{n}}{x_{n}}y_{n}v\ \d t
+\int_{\{|x_{n}(t)|<M\}}\frac{u_{n}}{\|x_{n}\|}v\ \d t\right|\\
\ & \leq \varepsilon
\|y_{n}\|_{\infty}\|v\|_{1}+\int_{\{|x_{n}(t)<M|\}}
\frac{a_{1}(t)(1+M^{r-1})}{\|x_{n}\|}|v|\d t\\
&\hskip 5cm (\t{see hypothesis $\hbox{H}(\hbox{j})_{5}(\hbox{iii})$})\\
\Rightarrow &\mathop {\lim \sup }\limits_{n \to \infty }\left|
\int_{0}^{b}\frac{u_{n}}{\|x_{n}\|}v\ \d t\right|
\leq\varepsilon\|y\|_{\infty}\|v\|_{1}.
\end{align*}

Since $\varepsilon >0$ was arbitrary, we infer that $\mathop
{\lim }_{n \to \infty}
\int_{0}^{b}\frac{u_{n}}{\|x_{n}\|}v \d t =0$. Also we have
$\mathop {\lim}_{n \to \infty }
\int_{0}^{b}\frac{h}{\|x_{n}\|}v \d t =0$. So from (23) and if
$v=y_{n}-y$, we obtain
\begin{align*}
\ &\mathop {\lim}\limits_{n \to \infty }
\left[\int_{0}^{b}y_{n}'(y_{n}-y)'\d t-\lambda_{m}\int_{0}^{b}y_{n}(y_{n}-y)\d t\right]=0,\\
\Rightarrow &\mathop {\lim}\limits_{n \to \infty }
\int_{0}^{b}y_{n}'(y_{n}'-y')\d t  \ \ (\t{recall that $y_{n}\rightarrow y$ in $C_{\rm per}(T)$}), \\
\Rightarrow & \;\|y_{n}'\|_{2}\rightarrow \|y'\|_{2}.
\end{align*}

Because $y_{n}'\xrightarrow{w} y'$ in $L^{2}(T)$, it follows that $y_{n}'\rightarrow y'$
in $L^{2}(T)$ and so $y_{n}\rightarrow y$ in $W_{\rm per}^{1,2}(T)$.

Since ${|\varphi
(x_{n})|}/{\|x_{n}\|^{2}}\leq {M_{1}}/{\|x_{n}\|^{2}}$ and
$\int_{0}^{b} ({j(t,x_{n}(t))}/{\|x_{n}\|^{2}})\d t\rightarrow 0$
(see hypothesis H(j)$_{5}$(v)), if in (23) $v=y_{n}$ and we pass
to the limit as $n\rightarrow\infty$, we obtain that
$\|y'\|_{2}^{2}= \lambda_{m}\|y\|_{2}^{2}$, hence
$y(t)=\xi_{1}\sin m\omega t+\xi_{2} \cos m\omega t$ with
$\xi_{1},\xi_{2}\in \mathbb{R}$ and so $y(t)=r\sin(m\omega
t+\theta)$ with $r=(\xi_{1}^{2}+ \xi_{2}^{2})^{\frac{1}{2}}$,
$\tan\theta=\xi_{1}/\xi_{2}$.

We write $y_{n}=\bar{y}_{n}+y_{n}^{0}+\hat{y}_{n}$ with
$\bar{y}_{n}\in\bar{H}$, ${y}_{n}^{0} \in N_{m}$,
$\hat{y}_{n}\in\hat{H}$, $n\geq 1$. Using
$v=-\bar{y}_{n}+y_{n}^{0}+\hat{y}_{n}\in W_{\rm per}^{1,2}(T)$ as our
test function, we obtain
\begin{align}
\ &
\;\;\left|\int_{0}^{b}x_{n}'(-\bar{y}_{n}+y_{n}^{0}+\hat{y}_{n})'\d t-
\lambda_{m}\int_{0}^{b}x_{n}(-\bar{y}_{n}+y_{n}^{0}+\hat{y}_{n})\d t \right. \nonumber\\
\ & \left.-\int_{0}^{b}u_{n}(-\bar{y}_{n}+y_{n}^{0}+\hat{y}_{n})\d t +
\int_{0}^{b}h(-\bar{y}_{n}+y_{n}^{0}+\hat{y}_{n})\d t \right|\nonumber \\
\leq &\; \varepsilon_{n}\|-\bar{y}_{n}+y_{n}^{0}+\hat{y}_{n}\|\leq
3\varepsilon_{n} \quad
\t{with} \ \ \varepsilon_{n}\downarrow 0 \nonumber\\
\Rightarrow &
\;\frac{1}{\|x_{n}\|}\left|\int_{0}^{b}x_{n}'(-\bar{x}_{n}'+x_{n}^{0'}+\hat{x}'_{n})\d t-
\lambda_{m}\int_{0}^{b}x_{n}(-\bar{x}_{n}+x_{n}^{0}+\hat{x}_{n})\d t\right|\nonumber \\
\ &-\left|\int_{0}^{b}\left(\frac{u_{n}-h}{\|x_{n}\|}\right)(-\bar{x}_{n}+x_{n}^{0}+\hat{x}_{n})\d t\right|
\leq 3 \varepsilon_{n}
\end{align}
where $x_{n}=\bar{x}_{n}+x_{n}^{0}+\hat{x}_{n}$ with
$\bar{x}_{n}\in\bar{H}$, $x_{n}^{0} \in N_{m}$,
$\hat{x}_{n}\in\hat{H}$. Because $x_{n}^{0}\in N_{m}$, we have
$x_{n}^{0}=\xi_{n} ^{1}\sin m\omega t + \xi_{n}^{2}\cos m\omega
t$ and so for all $n\geq 1$ we have $\|x_{n}^{0'}\|
_{2}^{2}=\lambda_{m}\|x_{n}^{0}\|_{2}^{2}$. Since
$|\langle x_{n}^{*},v\rangle |\leq\varepsilon_{n}\|v\|$ for all $v\in
W_{\rm per}^{1,2}(T)$, taking $v=x_{n}^{0}$ and exploiting the
orthogonality relations, we have
\begin{align}
\ &\;
\left|\frac{1}{\|x_{n}\|}(\|x_{n}^{0'}\|_{2}^{2}-\lambda_{m}\|x_{n}^{0}\|_{2}^{2})-\int_{0}^{b}
\frac{u_{n}}{\|x_{n}\|}x_{n}^{0}\d t+\int_{0}^{b}
\frac{h}{\|x_{n}\|}x_{n}^{0}\d t\right|<\varepsilon_{n},\nonumber\\
\Rightarrow &\;\int_{0}^{b}
\frac{u_{n}-h}{\|x_{n}\|}x_{n}^{0}\d t\rightarrow 0.
\end{align}

Also from Lemma~3 of [33], we have that
$\beta_{1}\|\bar{x}_{n}\|^{2}\leq\lambda_{m}
\|x_{n}\|_{2}^{2}-\|\bar{x}_{n}'\|_{2}^{2}$ for all $n\geq 1$
and some $\beta_{1}>0$ while from Lemma~9, we have
$\beta_{1}\|\hat{x}_{n}\|^{2}\leq\|\hat{x}_{n}'\|_{2}^{2}-
\lambda_{m}\|\hat{x}_{n}'\|^{2}_{2}$ for all $n\geq 1$. Using the
orthogonality relations among the three subspaces $\bar{H}$,
$N_{m}$, and $\hat{H}$, we obtain
\begin{align*}
\int_{0}^{b}x_{n}'(-\bar{x}_{n}'+x_{n}^{0'}+\hat{x}_{n}')\d t
&=-\|\bar{x}_{n}'\|_{2}^{2}+
\|x_{n}^{0'}\|_{2}^{2}+\|\hat{x}_{n}'\|_{2}^{2}\\
\intertext{and}
-\lambda_{m}\int_{0}^{b}x_{n}(-\bar{x}_{n}+x_{n}^{0}+\hat{x}_{n})\d t &=
\lambda_{m}\|\bar{x}_{n}\|^{2}_{2}-\lambda_{m}\|x_{n}^{0}\|_{2}^{2}-\lambda_{m}
\|\hat{x}_{n}\|_{2}^{2}.
\end{align*}

Thus finally we can write that
\begin{align}
\ &\int_{0}^{b}x_{n}'(-\bar{x}_{n}'+x_{n}^{0'}+\hat{x}_{n}')\d t-
\lambda_{m}\int_{0}^{b}x_{n}(-\bar{x}_{n}+x_{n}^{0}+\hat{x}_{n})\d t \nonumber \\[6pt]
=&\;\lambda_{m}\|\bar{x}_{n}\|^{2}_{2}-\|\bar{x}_{n}'\|^{2}_{2}+\|\hat{x}_{n}'\|_{2}^{2}
-\lambda_{m}\|\hat{x}_{n}\|_{2}^{2} \nonumber\\[6pt]
\geq &\; \beta_{2}\|-\bar{x}_{n}+\hat{x}_{n}\|^{2}\quad \t{for all}
\ \ n\geq 1 \ \ \t{and some} \ \ \beta_{2}>0.
\end{align}

Moreover, recalling that for almost all $t\in T$, all $|x|\geq M$
and all $u\in\partial j(t,x)$ we have
$|{u}/{x}|\leq\varepsilon$ and using also hypothesis
H(j)$_{5}$(iii), we have
\begin{align}
\int_{0}^{b}\frac{u_{n}}{\|x_{n}\|}(-\bar{x}_{n}+\hat{x}_{n})\d t &=
\frac{1}{\|x_{n}\|}\int_{\{|x_{n}(t)|\geq M\}}\frac{u_{n}}{x_{n}}x_{n}(-\bar{x}_{n}+\hat{x}_{n})\d t\nonumber\\
&\quad\,+\int_{\{|x_{n}(t)|< M\}}\frac{u_{n}}{\|x_{n}\|}(-\bar{x}_{n}+\hat{x}_{n})\d t \nonumber\\
\ & \leq\frac{1}{\|x_{n}\|}\varepsilon(\|\bar{x}_{n}\|_{2}^{2}+\|\hat{x}_{n}\|_{2}^{2})+
\frac{1}{\|x_{n}\|}\beta_{3}\|-\bar{x}_{n}+x_{n}\|\nonumber\\[3pt]
&\quad\,\t{for some} \ \ \beta_{3}>0 \nonumber\\
\ &\leq\frac{\varepsilon}{\|x_{n}\|}\|w_{n}\|^{2}+\frac{1}{\|x_{n}\|}\beta_{3}
\|w_{n}\|\nonumber\\[3pt]
&\,\quad \t{with} \ \ w_{n}=-\bar{x}_{n}+\hat{x}_{n}, n\geq 1.
\end{align}

Also for all $n\geq 1$, we have
\begin{align}
\left|\int_{0}^{b}h(-\bar{y}_{n}+\hat{y}_{n})\d t\right|&\leq\|h\|_{1} \ \|-\bar{y}_{n}+\hat{y}_{n}\|_{\infty}
\leq\beta_{4}\|-\bar{y}_{n}+\hat{y}_{n}\|\nonumber\\
&=\frac{\beta_{4}}{\|x_{n}\|}\|w_{n}\| \quad
\t{for some} \ \ \beta_{4}>0.
\end{align}

Using (25)$\to$(28) in (14), we obtain
\begin{equation*}
\hskip -2pc \frac{1}{\|x_{n}\|}(\beta_{2}-\varepsilon)\|w_{n}\|^{2}-\frac{1}{\|x_{n}\|}\beta_{5}\|w_{n}\|
\leq\varepsilon_{n}' \quad\t{for some} \ \ \beta_{5}>0 \ \ \t{and with} \ \
\varepsilon_{n}'\downarrow 0.
\end{equation*}

Choosing $\varepsilon <\beta_{2}$, we have
\begin{align*}
\ &\;\mathop {\lim \sup }\limits_{n \to \infty }
\frac{1}{\|x_{n}\|}(\beta_{6}\|w_{n}\|^{2}-\beta_{5}\|w_{n}\|)\leq
0 \quad \t{with} \ \
\beta_{6}=\beta_{2}-\varepsilon >0\\
\Rightarrow &\;\mathop {\lim \sup }\limits_{n \to \infty }
\frac{\|w_{n}\|^{2}}{\|x_{n}\|}\left(\beta_{6}-\frac{\beta_{5}}{\|w_{n}\|}\right)\leq 0,\\
\Rightarrow &\;\frac{\|w_{n}\|^{2}}{\|x_{n}\|}\rightarrow 0 \ \
(\t{by passing to a subsequence if necessary}).
\end{align*}

From the choice of the sequence $\{x_{n}\}_{n\geq 1}\subseteq W_{\rm per}^{1,2}(T)$ we have
\begin{equation*}
|\varphi(x_{n})|=\left|\frac{1}{2}\|x_{n}'\|_{2}^{2}-\frac{\lambda_{m}}{2}\|x_{n}\|_{2}^{2}-
\int_{0}^{b}j(t,x_{n}(t))\d t+\int_{0}^{b}h(t)x_{n}(t)\d t \right|\leq M_{1}.
\end{equation*}

Divide by $\|x_{n}\|$ and use the orthogonality of the subspaces and the equality
$\|x_{n}^{0'}\|_{2}^{2}=\lambda_{m}\|x_{n}^{0}\|_{2}^{2}$ for all $n\geq 1$. We obtain
\begin{align}
\frac{|\varphi(x_{n})|}{\|x_{n}\|}&=\left|\frac{1}{2}\frac{\|w_{n}'\|_{2}^{2}}{\|x_{n}\|}
-\frac{\lambda_{m}}{2}\frac{\|w_{n}\|_{2}^{2}}{\|x_{n}\|}-
\int_{0}^{b}\frac{j(t,x_{n}(t))}{\|x_{n}\|}\d t\right.\nonumber\\
&\quad\,\left.+\int_{0}^{b}h(t)y_{n}(t)\d t \right|\leq \frac{M_{1}}{\|x_{n}\|}.
\end{align}

Recall that $W^{1,2}_{\rm per}(T)=\mathbb{R}\oplus V$ with $V=\{v\in
W^{1,p}_{\rm per}(T)\!\!:\int_0^bv(t)\d t=0\}=\mathbb{R}^{\bot}$. So
$w_{n}=\xi_{n}+v_{n}$ with $\xi_{n}\in\mathbb{R}$, $v_{n}\in V$,
$n\geq 1$. We have $\|w_{n}\|_{2}^{2}=
\|\xi_{n}\|_{2}^{2}+\|v_{n}\|_{2}^{2}=b\xi_{n}^2+\|v_{n}\|_{2}^{2}$,
$n\geq 1$. So
\begin{align*}
&\;\frac{\lambda_{m}}{2}\|w_{n}\|_{2}^{2}=\frac{\lambda_{m}}{2}
\xi_{n}^{2}b+\frac{\lambda_{m}}{2}\|v_{n}\|_{2}^{2},\\[.2pc]
\Rightarrow
&\;\frac{\lambda_{m}}{2}\frac{\|w_{n}\|_{2}^{2}}{\|x_{n}\|}=
\frac{\lambda_{m}}{2}
\frac{\xi_{n}^{2}b}{\|x_{n}\|}+\frac{\lambda_{m}}{2}\frac{\|v_{n}\|_{2}^{2}}{\|x_{n}\|}
\leq\lambda_{m}\frac{\|w_{n}\|^{2}}{\|x_{n}\|}\rightarrow 0.
\end{align*}

Also we have
\begin{align*}
\
&\;\frac{1}{2}\|w_{n}'\|_{2}^{2}=\frac{1}{2}\|v_{n}'\|_{2}^{2}\leq\frac{\beta_{7}}{2}
\quad \t{with} \ \ \beta_{7}>0 \ \\[.2pc]
&\quad (\t{by the Poincar\'e--Wirtinger inequality}),\\[.2pc]
\Rightarrow &\;\frac{1}{2}\frac{\|w_{n}'\|_{2}^{2}}{\|x_{n}\|}
\leq\frac{\beta_{7}}{2}\frac{\|w_{n}\|^{2}}{\|x_{n}\|}\rightarrow
0 \quad \t{as} \ \ n\rightarrow\infty.
\end{align*}

Moreover, we have
\begin{equation*}
\int_{0}^{b}\frac{j(t,x_{n}(t))}{\|x_{n}\|}\d t=\int_{\{x_{n}(t)\neq 0\}}
\frac{j(t,x_{n}(t))}{x_{n}}y_{n}(t)\d t+\int_{\{x_{n}(t)= 0\}}
\frac{j(t,0)}{\|x_{n}\|}\d t.
\end{equation*}

Note that on $\{t\in T\!\!:y(t)>0\}$ we have $x_{n}(t)\rightarrow
+\infty$ and on $\{t\in T\!\!:y(t)<0\}$ we have that
$x_{n}(t)\rightarrow -\infty$. In addition $\int_{\{x_{n}(t)= 0\}}
\frac{j(t,0)}{\{\|x_{n}\|\}}\d t\rightarrow 0$. So via Fatou's
lemma, we have
\begin{align*}
\mathop {\lim \inf }\limits_{n \to \infty }\int_{0}^{b}\frac{j(t,x_{n}(t))}{\|x_{n}\|}\d t &
\geq\mathop {\lim \inf }\limits_{n \to \infty }\int_{\{x_{n}(t)\neq 0\}}
\frac{j(t,x_{n}(t))}{x_{n}(t)}y_{n}(t)\d t \\[.2pc]
\ & \geq\int_{0}^{b}j_{+}(t)y^{+}(t)\d t -\int_{0}^{b}j_{-}(t)y^{-}(t)\d t.
\end{align*}

Because $\frac{\|w_{n}\|}{\|x_{n}\|}\rightarrow 0$, we have that $y\in N_{m}$ and so
\begin{align}
\mathop {\lim \inf }\limits_{n \to \infty }\int_{0}^{b}\frac{j(t,x_{n}(t))}{\|x_{n}\|}\d t &
\geq\int_{0}^{b}j_{+}(t)r\sin(m\omega t+\theta)^{+}\d t \nonumber\\[.2pc]
\ &\;\;\;\;-\int_{0}^{b}j_{-}(t)r\sin(m\omega t+\theta)^{-}\d t,
\theta\in\mathbb{R}.
\end{align}

From (29) and since $\frac{\|w_{n}\|^{2}}{\|x_{n}\|}\rightarrow
0$, we obtain
\begin{align}
&\mathop {\lim }\limits_{n \to \infty }\int_{0}^{b}\frac{j(t,x_{n}(t))}{\|x_{n}\|}\d t
=\;\int_{0}^{b}h(t)y(t)\d t =\int_{0}^{b}h(t)r\sin(m\omega t+\theta)\d t \nonumber\\[.2pc]
\ &\quad <\;r \bigg(\int_{0}^{b}j_{+}(t)\sin(m\omega t+\theta)^{+}\d t
-\int_{0}^{b}j_{-}(t)\sin(m\omega t+\theta)^{-}\d t \bigg)\nonumber\\[.2pc]
\ &\quad\;\;\;\;\;(\t{by hypothesis} \ \ \hbox{H}(\hbox{j})_{5}(\hbox{v})).
\end{align}

Comparing (30) and (31), we reach a contradiction. This proves
that the sequence $\{x_{n}\}_{n\geq 1}\subseteq W_{\rm per}^{1,2}(T)$
is bounded, hence we may assume that $x_{n}\xrightarrow{w} x$ in
$W_{\rm per}^{1,2}(T)$ and $x_{n}\rightarrow x$ in $C_{\rm per}(T)$. As
before (see the proof of Theorem 1) we can finish the proof and
conclude that $\varphi$ satisfies the non-smooth PS-condition.\hfill QED
\end{proof}

Let $H_{1}=\t{span}\{\sin k\omega t,\cos k\omega t\!\!:k=0,1,...,m\}$
and $H_{2}=H_{1}^{\bot}$.

\begin{prop}$\left.\right.$\vspace{.5pc}

\noindent If hypotheses ${\rm H(j)}_{5}$ hold{\rm ,}
then $\varphi(x)\rightarrow -\infty$ as $\|x\|\rightarrow \infty${\rm ,} $x\in H_{1}$.
\end{prop}

\begin{proof}
Suppose that the conclusion of the proposition was not true. Then we
can find $\beta\in\mathbb{R}$ and a sequence $\{x_{n}\}_{n\geq 1}\subseteq H_{1}$ such that
$\|x_{n}\|\rightarrow\infty$ and $\varphi(x_{n})\geq\beta$ and all $n\geq 1$. We have
\begin{equation*}
\frac{1}{2}\|x_{n}'\|_{2}^{2}-\frac{\lambda_{m}}{2}\|x_{n}\|_{2}^{2}-
\int_{0}^{b}j(t,x_{n}(t))\d t+\int_{0}^{b}h(t)x_{n}(t)\d t\geq\beta.
\end{equation*}

Let $y_{n}=\frac{x_{n}}{\|x_{n}\|}$, $n\geq 1$. We may assume that $y_{n}\xrightarrow{w} y$
in $W_{\rm per}^{1,2}(T)$ and $y_{n}\rightarrow y$ in $C_{\rm per}(T)$. Because $H_{1}$ is finite
dimensional and $y_{n}\in H_{1}$ for all $n\geq 1$, we have $y_{n}\rightarrow y$ in
$W_{\rm per}^{1,2}(T)$. For all $n\geq 1$ we have
\begin{align*}
\hskip .05pc &\frac{1}{2}\|y_{n}'\|_{2}^{2}-\frac{\lambda_{m}}{2}\|y_{n}\|_{2}^{2}-
\int_{0}^{b}\frac{j(t,x_{n}(t))}{\|x_{n}\|^{2}}\d t+\int_{0}^{b}\frac{h(t)}{\|x_{n}\|}
y_{n}(t)\d t\geq\frac{\beta}{\|x_{n}\|^{2}},\\
&\quad\,n\geq 1.
\end{align*}

Clearly $\int_{0}^{b}\frac{j(t,x_{n}(t))}{\|x_{n}\|^{2}}\d t\rightarrow 0$ and
$\int_{0}^{b}\frac{h(t)}{\|x_{n}\|}y_{n}(t)\d t\rightarrow 0$. So in the limit as
$n\rightarrow\infty$, we obtain
\begin{align*}
\ &\;\frac{\lambda_{m}}{2}\|y\|_{2}^{2}\leq\frac{1}{2}\|y'\|_{2}^{2}\\
\Rightarrow
&\;\frac{1}{2}\|y'\|_{2}^{2}=\frac{\lambda_{m}}{2}\|y_{n}\|_{2}^{2}
\ \
(\t{since $y\in H_{1}$})\\
\Rightarrow &\;y\in N_{m}.
\end{align*}

From the choice of the sequence $\{x_{n}\}_{n\geq 1}\subseteq H_{1}$ we have
\begin{align*}
\ &-\int_{0}^{b}j(t,x_{n}(t))\d t+\int_{0}^{b}h(t)x_{n}(t)\d t\geq\varphi(x_{n})\geq\beta\\
\Rightarrow &\int_{0}^{b}\frac{j(t,x_{n}(t))}{\|x_{n}\|}\d t\leq-\frac{\beta}{\|x_{n}\|}+
\int_{0}^{b}h(t)y_{n}(t)\d t.
\end{align*}

Arguing as in the proof of Proposition~10, in the limit as
$n\rightarrow \infty$, we obtain
\begin{align*}
\ &\;r\int_{0}^{b}j_{+}(t)\sin(m\omega t+\theta)^{+}\d t
-r\int_{0}^{b}j_{-}(t)\sin(m\omega t+\theta)^{-}\d t\\
\leq &\;r\int_{0}^{b}h(t)\sin(m\omega t+\theta)\d t,
\end{align*}
which contradicts hypothesis H(j)$_{5}$(v). This proves the
proposition.\hfill QED
\end{proof}

\begin{prop}$\left.\right.$\vspace{.5pc}

\noindent If hypotheses ${\rm H(j)}_{5}$ hold{\rm ,}
then $\varphi(x)\rightarrow +\infty$ as $\|x\|\rightarrow \infty${\rm ,} $x\in H_{2}$.
\end{prop}

\begin{proof}
For $x\in H_{2}$, we have
\begin{align*}
\varphi(x) &=\frac{1}{2}\|x'\|_{2}^{2}-\frac{\lambda_{m}}{2}\|x\|_{2}^{2}
-\int_{0}^{b}j(t,x(t))\d t+\int_{0}^{b}h(t)x(t)\d t\\
\ &\geq c \|x\|^{2}-\int_{0}^{b}j(t,x(t))\d t+\int_{0}^{b}h(t)x(t)\d t  \ \ (\t{see Lemma 9})\\
\Rightarrow\frac{\varphi(x) }{\|x\|^{2}} &\geq c-\int_{0}^{b}\frac{j(t,x(t))}{\|x\|^{2}}\d t +
\int_{0}^{b}\frac{h(t)}{\|x\|}y(t)\d t.
\end{align*}

Remark that $\int_{0}^{b}\frac{j(t,x(t))}{\|x\|^{2}}\d t\rightarrow 0$ and
$\int_{0}^{b}\frac{h(t)}{\|x\|}y(t)\d t\rightarrow 0$ as $\|x\|\rightarrow\infty$, $x\in H_{2}$.
So we have
\begin{align*}
\ &\;\mathop {\lim \inf }\limits_{\mathop {\|x\| \to \infty
}\limits_{x \in H_{ 2} } }
\frac{\varphi(x) }{\|x\|^{2}} \geq c>0 \\
\Rightarrow &\;\varphi(x)\rightarrow +\infty \quad \t{as} \ \
\|x\|\rightarrow\infty,x\in H_{2}.
\end{align*}
$\left.\right.$\hfill QED
\end{proof}

Propositions~10--12, permit the use of the non-smooth saddle
point theorem which gives $x\in W_{\rm per}^{1,2}(T)$ such that
$0\in\partial\varphi(x)$. As in the proof of Theorem~1, we can
check that $x\in C^{1}(T)$, $x'\in W^{1,1}(T)$ and also $x$
solves (22). So we can state the following existence theorem.

\setcounter{lem}{12}
\begin{thm}[\!]
If hypotheses ${\rm H(j)}_{5}$ hold{\rm ,}
then for every $h\in L^{1}(T)${\rm ,} problem {\rm (22)} has
a solution $x\in C^{1}(T)$ with $x'\in W^{1,1}(T)$.
\end{thm}

\begin{rem}
{\rm Theorem~13 generalizes Theorems~2 and
3 of Tang [33]. The generalization is two-fold. On the one hand we
assume a more general Landesmann--Lazer type condition (see
hypothesis H(j)$_{5}$(v) and Proposition 8) and on the other hand
we have a non-smooth potential function. Moreover, in Tang the
potential function is independent of the time-variable $t\in T$.}
\end{rem}

\end{document}